\documentclass[preprint]{elsarticle}

\usepackage{lineno,hyperref}
\usepackage{amsmath, amssymb,epsfig}
\usepackage{enumerate} % added by GN
\usepackage{multirow} % added by GN
\usepackage{placeins}
\modulolinenumbers[5]

\journal{Journal of Computational and Applied Mathematics}

\newtheorem{thm}{Theorem}
\newtheorem{lem}{Lemma}
\newtheorem{prop}{Proposition}
\newdefinition{rmk}{Remark}
\newdefinition{dfn}{Definition}
\newdefinition{ex}{Example}
\newproof{pf}{Proof}
\newcommand{\ds}{\displaystyle} % added by GN
\begin{document}

\begin{frontmatter}

\title{Asymptotically optimal definite quadrature formulae of 4-th
order\tnoteref{This work was supported by the Sofia University
Research Fund through Contract no. 30/2016. }}
\tnotetext[mytitlenote]{This work was supported by the Sofia
University Research Fund through Contract no. 30/2016.}


\author{Ana Avdzhieva}
\ead{avdzhieva@fmi.uni-sofia.bg}
%[mymainaddress]
\author{Geno Nikolov\corref{mycorrespondingauthor}}
\cortext[mycorrespondingauthor]{Corresponding author}
\ead{geno@fmi.uni-sofia.bg}

\address{Faculty of Mathematics and Informatics,
         Sofia University ``St. Kliment Ohridski'',
         5~James Bourchier Blvd., 1164 Sofia, Bulgaria}

\begin{abstract}
We construct several sequences of asymptotically optimal definite
quadrature formulae of fourth order and evaluate their error
constants. Besides the asymptotical optimality, an advantage of our
quadrature formulae is the explicit form of their weights and nodes.
For the remainders of our quadrature formulae monotonicity
properties are established when the integrand is a 4-convex
function, and a-posteriori error estimates are proven.
\end{abstract}


\end{frontmatter}

%\linenumbers

\section{Introduction}

We study quadrature formulae of the form
\begin{equation}\label{e1}
Q_n[f]=\sum_{i=1}^{n}a_{i,n}\,f(\tau_{i,n})\,,\quad 0\le
\tau_{1,n}<\tau_{2,n}<\cdots<\tau_{n,n}\leq 1
\end{equation}
for approximate evaluation of the definite integral
$$
I[f]:=\int\limits_{0}^{1}f(x)\,dx\,.
$$
Our interest is in definite quadrature formulae. Let us recall some
definitions.
\begin{dfn}
Quadrature formula \eqref{e1} is said to be \emph{definite of order
$r$}, $\,r\in \mathbb{N}$, if there exists a real non-zero constant
$c_{r}(Q_n)$ such that its remainder functional admits the
representation
$$
R[Q_n;f]:=I[f]-Q_n[f]=c_r(Q_n)\,f^{(r)}(\xi)
$$
for every $f\in C^{r}[0,1]$, with some $\xi\in [0,1]$ depending on
$f$\,.\smallskip

Furthermore, $Q_n$ is called positive definite (resp., negative
definite) of order $r$, if $c_r(Q_n)>0$ ($c_r(Q_n)<0$).
\end{dfn}

Obviously, if $Q_n$ is a definite quadrature formula of order $r$,
then $Q_n$ has \emph{algebraic degree of precision} $r-1$ (in short,
$ADP(Q_n)=r-1$), i.e., $R[Q_n;f]=0$ whenever $f$ is an algebraic
polynomial of degree at most $r-1$, and $R[Q_n;x^{r}]\ne
0$.\smallskip

Throughout this paper, by $r$-convex ($r$-concave) function $f$ we
shall mean a function  $f\in C^{r}[0,1]$ such that $f^{(r)}\geq 0$
($f^{(r)}\leq 0$) on the interval $[0,1]$.\smallskip

The importance of definite quadrature formulae of order $r$ lies in
the one-sided approximation they provide for $I[f]$ when the
integrand $f$ is $r$-convex (concave). If, e.g., $\{Q^{+},\ Q^{-}\}$
is a pair of a positive and a negative definite quadrature formula
of order $r$ and $\,f\,$ is $r$-convex, then for the true value of
$I[f]$ we have the inclusion $Q^{+}[f]\leq I[f]\leq Q^{-}[f]$. This
simple observation serves as a base for derivation of a posteriori
error estimates and rules for termination of calculations (stopping
rules) in automatic numerical integration algorithms (see
\cite{KJF:1993} for a survey). Most of quadratures used in practice
(e.g., quadrature formulae of Gauss, Radau, Lobatto, Newton-Cotes)
are definite of certain order.\smallskip

Definite $n$-point quadrature formulae with smallest positive or
largest negative error constant are called optimal definite
quadrature formulae. Let us set
\begin{eqnarray*}
&&c_{n,r}^{+}:=\inf\{c_{n,r}(Q_n)\,:\, Q_n\ \text{ is positive
definite of order }\ r\}\,,\\
&&c_{n,r}^{-}:=\inf\{c_{n,r}(Q_n)\,:\, Q_n\ \text{ is negative
definite of order }\ r\}\,.
\end{eqnarray*}

It should be pointed out that it is fairly not obvious that the
above infimums are attained or that the optimal definite quadrature
formulae are unique. The existence of optimal definite quadrature
formulae was first proven by Schmeisser \cite{GS:1972} for even $r$,
and for arbitrary $r$ and more general boundary conditions by Jetter
\cite{KJ:1976} and Lange \cite{GL:1977}. The uniqueness has been
proven by Lange \cite{GL:1977, GL:1979}. For even $r$, Lange
\cite{GL:1977} has shown that
\begin{equation}\label{e2}
\begin{split}
&c_{n,r}^{+}=-\frac{B_r(j/2)}{n^{r}}\,\Big(1+O(n^{-1})\Big)\quad
\text{ if }\ r=4m+2j\,,\\
&c_{n,r}^{-}=-\frac{B_r(j/2)}{n^{r}}\,\Big(1+O(n^{-1})\Big)\quad
\text{ if }\ r=4m+2-2j
\end{split}
\end{equation}
for $j=1,\,2$, where $B_r$ is the $r$-th Bernoulli polynomial with
leading coefficient $1/r!$\,. Schmeisser \cite{GS:1972} proved that
the same result holds for optimal definite quadrature formulae with
equidistant nodes.\smallskip

The $n$-point optimal positive definite and the $(n+1)$-point
optimal negative definite quadrature formulae of order $2$ are
well-known: these are the $n$-th compound midpoint and trapezium
quadrature formulae, respectively. The case $r=2$ is exceptional, as
for $r\geq 3$ the optimal definite quadrature formulae are not
known. Lange \cite{GL:1977} has computed numerically, for $3\leq
n\leq 30$, the $n$-point optimal  definite quadrature formulae of
order $3$ and the $n$-point optimal positive definite quadrature
formulae of order $4$.\smallskip

It is a general observation about the optimality concept in
quadratures that, even though the existence and the uniqueness of
the optimal quadrature formulae (for instance, in the non-periodic
Sobolev classes of functions) is established, the optimal quadrature
formulae remain unknown. This fact severely reduces the practical
importance of optimal quadratures. The way out of this situation is
to look for quadrature formulae which are nearly optimal, e.g., for
sequences of asymptotically optimal quadrature formulae.

\begin{dfn}
Let $\,\{Q_n\}_{n=n_{0}}^{\infty}$ be a sequence of positive (resp,
negative) definite quadrature formulae of order $\,r$. $\,Q_n\,$ is
said to be \emph{asymptotically optimal} positive (negative)
definite quadrature formula of order $r$, if
$$
\lim_{n\rightarrow\infty}\frac{c_r(Q_n)}{c_{n,r}^{+}}=1\qquad
\text{resp., }\ \
\lim_{n\rightarrow\infty}\frac{c_r(Q_n)}{c_{n,r}^{-}}=1\,.
$$
\end{dfn}

In \cite{GS:1972} Schmeisser proposed an approach for construction
of asymptotically optimal definite quadrature formulae of even order
$r$ with equidistant nodes. K\"{o}hler and Nikolov \cite{KN:1995a}
have studied Gauss-type quadrature formulae associated with spaces
of splines with double and equidistant knots, and as a result
obtained bounds for the best constants $c_{n,r}^{+}$ and
$c_{n,r}^{-}$. In particular, it has been shown in \cite{KN:1995a}
that for even $r$ the corresponding Gauss-type quadrature formulae
are asymptotically optimal definite quadrature formulae. Motivated
by this result, in \cite{GN:1996} Nikolov found explicit recurrence
formulae for the evaluation of the nodes and the weights of the
Gaussian formulae for the spaces of cubic splines with double
equidistant knots, and proposed a numerical procedure for the
construction of the Lobatto quadrature formulae for the same spaces
of splines. According to \cite{KN:1995a}, the Gauss and the Lobatto
quadrature formulae for these spaces of splines are respectively
asymptotically optimal positive definite and asymptotically optimal
negative definite, of order $4$.\smallskip

Although the evaluation of Gauss-type quadrature formulae for spaces
of splines (also with single knots, because of their asymptotical
optimality in certain Sobolev classes, see \cite{KN:1995}) is highly
desirable, there is a serious problem occurring already with the
splines of degree $3$, and its difficulty increases with the splines
degree: the mutual displacement of the nodes of the quadratures and
the splines knots is unknown. For justifying the location of the
quadrature abscissae with respect to the knots of the space of
splines, additional assumptions are to be made. For instance, in a
recent paper \cite{ABC:2015}  Ait-Haddou, Barto\v{n} and Calo
extended the procedure from \cite{GN:1996} for explicit evaluation
of the Gaussian quadrature formulae for spaces of $C^{1}$ cubic
splines with non-equidistant knots, assumming that the spline knots
are symmetrically stretched. For another approach to the
construction of Gaussian quadrature formulae for $C^{2}$ cubic
splines via homotopy continuation, see \cite{BC:2016}.\smallskip

In the present paper we construct several sequences of
asymptotically optimal definite quadrature formulae of order $4$
with explicitly given nodes and weights. For their construction we
make use of the Euler-Maclaurin summation formulae, associated with
the midpoint and the trapezium quadrature formulae, replacing the
values of the derivatives at the end-points by appropriate formulae
for numerical differentiation (this idea is not new, it can be
traced in the book of Brass \cite{HB:1977}, and, implicitly, has
been applied already in \cite{GS:1972}). Thus, our quadrature
formulae differ from the compound midpoint or compound trapezium
quadrature formula by very little: they have only few different
weights and/or involve few additional nodes. We evaluate the error
constants of our quadrature formulae, which, in view of their
asymptotical optimality, are not essentially different. Our
motivation for proposing not just two sequences of asymptotically
optimal positive definite and negative definite quadrature formulae
of order $4$ is that, when chosen appropriately, pairs of definite
quadrature formulae of the same type furnish, similarly to the case
of pairs of definite quadrature formulae of opposite type, error
inclusions for $I[f]$ whenever the integrand $f$ is $4$-convex or
$4$-concave.\smallskip

The rest of the paper is organized as follows. Section~2 provides
the necessary facts about Peano representation theorem for linear
functionals, Bernoulli polynomials and Euler-Maclaurin summation
formulae. In Section~3 we construct sequences of definite quadrature
formulae of order $4$. In Section~4 we prove monotonicity of the
remainders of some of our definite quadrature formulae under the
assumption that the integrand is $4$-convex (concave), and as a
result obtain a posteriori error estimates. Section 5 shows some
numerical experiments, and Section 6 contains some final remarks.
%^================================================================
\section{Preliminaries}
\subsection{Peano kernel representation of linear functionals}

Throughout the paper, $\pi_{m}$ will stand for the set of algebraic
polynomials of degree not exceeding $m$.

By $W^{r}_1=W^{r}_1[0,1]$, $\,r\in \mathbb{N}$, we denote the
Sobolev class of functions
$$
W^{r}_1[0,1]:=\{f\in C^{r-1}[0,1]\,:\, f^{(r-1)} \mbox{ abs.
continuous},\; \int_{0}^{1}\!|f^{(r)}(t)|\,dt<\infty\}\,.
$$
In particular, we have $C^{r}[0,1]\subset W^{r}_1[0,1]$.

If ${\cal L}$ is a linear functional defined in $W^{r}_1[0,1]$ which
vanishes on $\pi_{r-1}$, then, by a classical result of Peano
\cite{GP:1913}, ${\cal L}$ admits the integral representation
$$
{\cal L}[f]=\int_{0}^{1}K_r(t)f^{(r)}(t)\,dt,\qquad \quad
K_r(t)={\cal L}\Big[\frac{(\cdot-t)_{+}^{r-1}}{(r-1)!}\Big],\ \ t\in
[0,1]\,,
$$
where
$$
u_{+}=\max\{u,0\}\,.
$$

The function $K_r$ is called \emph{the $r$-th Peano kernel of}
${\cal L}$. In the case when ${\cal L}$ is the remainder
$R[Q_n;\cdot]\,$ of the quadrature formula \eqref{e1} and
$ADP(Q_n)\geq r-1$, $\,K_r(t)=K_r(Q_n;t)\,$ is also referred to as
\emph{the $r$-th Peano kernel of $Q_n$}. An explicit representations
of $\,K_r(Q_n;t)\,$ for $t\in [0,1]$ is
\begin{equation}\label{e3}
K_r(Q_n;t)=(-1)^{r}\Big[\frac{t^r}{r!}-
\frac{1}{(r-1)!}\sum_{i=1}^{n}a_{i,n}(t-\tau_{i,n})_{+}^{r-1}\Big].
\end{equation}
$K_r(Q_,;\cdot)$ is called also \emph{a monospline} of degree $r$.
From
\begin{equation}\label{e4}
R[Q_n;f]=\int\limits_{0}^{1}K_r(Q_n;t)\,f^{(r)}(t)\,dt
\end{equation}
it is clear that $Q_n$ is a positive (negative) definite quadrature
formula of order $r$ if and only if $ADP(Q_n)=r-1$ and
$K_r(Q_n;t)\geq 0$ (resp., $K_r(Q_n;t)\leq 0$) on $[0,1]$. For more
details on the Peano kernel theory we refer to \cite{HB:1977}.

\subsection{Bernoulli polynomials. Summation formulae
of Euler--Maclaurin type}

By appropriate integration by parts in \eqref{e4} the remainder of a
quadrature formula $Q_n$ with $ADP(Q_n)=r-1\,$ can be further
expanded in the form
\begin{equation}\label{e5}
R[Q_n;f]=\sum_{\nu=r}^{s}C_{\nu}(0)\,\big[f^{(\nu-1)}(1)-f^{(\nu-1)}(0)\big]
+\int\limits_{0}^{1}C_s(t)\,f^{(s)}(t)\,dt\,,
\end{equation}
with some functions $\{C_{\nu}(x)\}$ depending on $Q_n$ (see
\cite{HB:1977} for details).

For the sake of convenience, let us fix some notations. For $n\in
\mathbb{N}$, we set
\begin{equation}\label{e6}
x_{k,n}:=\frac{k}{n}\,,\quad (0\leq k\leq n)\,,\qquad
y_{\ell,n}:=\frac{2\ell-1}{2n}\,,\quad (1\leq\ell\leq n)\,.
\end{equation}
The $n$-th compound trapezium and midpoint quadrature formulae are
denoted by $Q_n^{Tr}$ and $Q_n^{Mi}$, respectively, i.e.,
$$
Q_n^{Tr}[f]:=\frac{1}{2n}\big[f(x_{0,n})+f(x_{n,n})\big]+\frac{1}{n}\,
\sum_{k=1}^{n-1}f(x_{k,n})\,,\qquad
Q_n^{Mi}[f]:=\frac{1}{n}\,\sum_{k=1}^{n}f(y_{k,n})\,.
$$

The Bernoulli polynomials $B_{\nu}$ are defined recursively by
$$
B_0(x)=1,\quad B_{\nu}^{\prime}(x)=B_{\nu-1}(x),\quad \
\int_0^1B_{\nu}(t)\,dt=0\,, \qquad  \nu\in \mathbb{N}.
$$
Here, we shall need the explicit form of $B_4(x)$,
$B_4(x)=x^2(1-x)^2/24-1/720$, and shall exploit the fact that
\begin{equation}\label{e7}
-\frac{1}{720}=B_4(0)\leq B_4(x)\leq
B_4(1/2)=\frac{7}{5760}\,,\qquad x\in [0,1]\,.
\end{equation}

The $1$-periodic extension of $B_{\nu}(x)$ on $\mathbb{R}$ is
denoted by $\widetilde{B}_{\nu}(x)$ and is called \emph{Bernoulli
monospline}. The expansion \eqref{e5} with $\,Q_n=Q_n^{Mi}\,$ and
$\,Q_n=Q_n^{Tr}\,$ yields the so-called Euler-Maclaurin summation
formulae (see, e.g., \cite[Satz 98, 99]{HB:1977}). For easier
further reference, they are given in a lemma:
\begin{lem}\label{l1}
Assume that $f\in W^{s}_1$, where $s\in \mathbb{N},\ s\geq 2$. Then
$$
R[Q_n^{Mi};f]\!=\!-\!\sum_{\nu=1}^{[\frac{s}{2}]}\!
\frac{B_{2\nu}(1/2)}{n^{2\nu}}\big[f^{(2\nu\!-\!1)}(1)\!-\!f^{(2\nu\!-\!1)}(0)\big]
\!+\!\frac{(-1)^{s}}{n^{s}}\!\int\limits_0^{1}\!\!
\widetilde{B}_s\Big(nx-\frac{1}{2}\Big)f^{(s)}(x)dx
$$
and
$$
R[Q_{n}^{Tr};f]=-\!\sum_{\nu=1}^{[\frac{s}{2}]}
\frac{B_{2\nu}(0)}{n^{2\nu}}\big[f^{(2\nu\!-\!1)}(1)\!-\!f^{(2\nu\!-\!1)}(0)\big]
\!+\!\frac{(-1)^{s}}{n^{s}}\!\int\limits_0^{1}\!\widetilde{B}_s(nx)f^{(s)}(x)dx\,.
$$
\end{lem}
%==========================================================
\section{Asymptotically optimal definite quadrature formulae of $4$-th order}
For verifying the asymptotical optimality of the definite quadrature
formulae constructed in this section, we note that \eqref{e2} with
$r=4$ reads as
\begin{equation}\label{e8}
c_{n,4}^{+}=\frac{1}{720\,n^{4}}\,\Big(1+O(n^{-1})\Big)\,,\qquad
c_{n,4}^{-}=-\frac{7}{5760\,n^{4}}\,\Big(1+O(n^{-1})\Big)\,.
\end{equation}

\begin{dfn}\label{d3}
For a given $\mathbf{t}=(t_1,t_2,t_3,t_4)$, $0\leq
t_1<t_2<t_3<t_4<1/2$, we denote by $\,D_1(\mathbf{t})[f]\,$ and
$\,D_3(\mathbf{t})[f]\,$ the interpolatory formulae for numerical
differentiation with nodes $\{t_i\}_{i=1}^{4}$, which approximate
$f^{\prime}(0)$ and $f^{\prime\prime\prime}(0)$, respectively. The
formulae approximating $f^{\prime}(1)$ and
$f^{\prime\prime\prime}(1)$ and obtained by reflection are denoted
by $\,\widetilde{D}_1(\mathbf{t})[f]\,$ and
$\,\widetilde{D}_3(\mathbf{t})[f]\,$, i.e.,
$$
\widetilde{D}_k(\mathbf{t})[f]=D_k(\mathbf{t})[g]\,,\qquad
g(x)=-f(1-x)\,,\quad k=1,3\,.
$$
For the sake of brevity, we write $D(\mathbf{t})$ for the collection
of four formulae for numerical differentiation
$\{D_1[\mathbf{t}],\,\widetilde{D}_1(\mathbf{t}),\,
D_3[\mathbf{t}],\,\widetilde{D}_3(\mathbf{t})\}$\,.
\end{dfn}
%==============================================================
\subsection{Negative definite quadrature formulae of order $4$ based on $Q_n^{Tr}$}
The second formula in Lemma \ref{l1} with $s=4$ can be rewritten in
the form
\[
\begin{split}
\int\limits_{0}^{1}f(x)\,dx=& Q_{n}^{Tr}[f]
-\frac{1}{12\,n^2}\big[f^{\prime}(1)-f^{\prime}(0)\big]+
\frac{1}{384\,n^4}\big[f^{\prime\prime\prime}(1)-f^{\prime\prime\prime}(0)\big]\\
&+\frac{1}{n^{4}}\int\limits_0^{1}\big[\widetilde{B}_4(nx)-B_4(1/2)\big]\,f^{(4)}(x)\,dx\\
=:& Q_n^{\prime}[f]+
\frac{1}{n^{4}}\int\limits_0^{1}\big[\widetilde{B}_4(nx)-B_4(1/2)\big]\,f^{(4)}(x)\,dx\,.
\end{split}
\]
Clearly, $Q_n^{\prime}$ is a negative definite quadrature formula of
order $4$, since, in view of \eqref{e7},
$K_4(Q_n^{\prime};x)=n^{-4}\big[\widetilde{B}_4(nx)-B_4(1/2)\big]\leq
0$. However, $Q_n^{\prime}$ is not of the desired form, as it
involves derivatives of the integrand. Therefore, we choose a set
$D(\mathbf{t})$ of formulae for numerical differentiation to replace
the values of $f^{\prime}$ and $f^{\prime\prime\prime}$ in
$Q_n^{\prime}$, and thus to obtain a (symmetric) quadrature formula
\begin{equation}\label{e9}
Q=Q_n^{Tr}+\frac{1}{12\,n^2}\big(D_1(\mathbf{t})-
\widetilde{D}_1(\mathbf{t})\big)-\frac{1}{384\,n^{4}}\big(D_3(\mathbf{t})
-\widetilde{D}_3(\mathbf{t})\big)\,,
\end{equation}
which involves at most $8$ nodes in addition to
$\{x_{k,n}\}_{k=0}^{n}$. (In the sequel, we shall refer to $Q$ to as
a \emph{quadrature formula generated by $D(\mathbf{t})$}.) We have
$$
R[Q;f]=R[Q_n^{\prime};f]+\frac{1}{12\,n^2}\,(L_1[f]-\widetilde{L}_1[f])-
\frac{1}{384\,n^{4}}\,(L_3[f]-\widetilde{L}_3[f])\,,
$$
where
$$
L_k[f]:=f^{(k)}(0)-D_k(\mathbf{t})[f]\,,\quad
\widetilde{L}_k[f]:=f^{(k)}(1)-\widetilde{D}_k(\mathbf{t})[f]\,,\qquad
k=1,\,3\,.
$$
The linear functionals $L_k$ and $\widetilde{L}_k$ vanish on
$\pi_3$, hence $R[Q;f]$ vanishes on $f\in\pi_3$, too. From the
definition of the Peano kernels it is readily seen that
$$
K_4(L_k;t)\equiv 0,\ \ t\in (t_4,1]\,,\qquad
K_4(\widetilde{L}_k;t)\equiv 0,\ \ t\in [0,1-t_4)\,, \quad
k=1,\,3\,.
$$
This implies the following important observation:
\begin{prop}\label{p1}
The fourth Peano kernel of the symmetric quadrature formula $Q$
generated by $D(\mathbf{t})$ through \eqref{e9} satisfies
$$
K_4(Q;x)\equiv
\frac{1}{n^{4}}\,\big[\widetilde{B}_4(nx)-B_4(1/2)\big],\qquad x\in
[t_4,1-t_4]\,.
$$
As a consequence, $\,Q\,$ is negative definite of order $4$ if and
only if $\,K_4(Q;x)\leq 0\,$ for $\,x\in (0,t_4)$\,.
\end{prop}

It should be pointed out that not every set $D(\mathbf{t})$ of
formulae for numerical differentiation generates a definite
quadrature formula of order $4$.

Our first application of the above approach reveals a known result.
%=================================================================
\subsubsection{A quadrature formula of G. Schmeisser}
The following $(n+1)$-point (with $n\geq 7$) asymptotically optimal
negative definite quadrature formula of order $4$ was obtained in
\cite[eqn. (43)]{GS:1972}: {\small
\[
\begin{split}
Q_{n+1}[f]=&\frac{403}{1152n}\big[f(x_{0,n})\!+\!f(x_{n,n})\big]\!+\!
\frac{159}{128n}\big[f(x_{1,n})\!+\!f(x_{n-1,n})\big]\\
&\!+\!\frac{113}{128n}\big[f(x_{2,n})\!+\!f(x_{n-2,n})\big]\!+\!
\frac{1181}{1152n}\big[f(x_{3,n})\!+\!f(x_{n-3,n})\big]
\!+\!\frac{1}{n}\sum_{k=4}^{n-4}f(x_{k,n})\,,
\end{split}
\]}
with an error constant
$$
c_4(Q_{n+1})=-\frac{7}{5760\,n^{4}}\Big(1+\frac{195}{7n}\Big)\,.
$$
Schmeisser's quadrature formula is generated by
$D(x_{0,n},x_{1,n},x_{2,n},x_{3,n})$. As this is a known result, we
do not enter into details.
%=============================================================
\subsubsection{A quadrature formula generated by
$D(x_{0,n},x_{1,3n},x_{2,3n},x_{1,n})$}

For $\mathbf{t}=(x_{0,n},x_{1,3n},x_{2,3n},x_{1,n})$ we have
$$
D_1(\mathbf{t})[f]=\frac{n}{2}\big[-11f(x_{0,n})
+18f(x_{1,3n})-9f(x_{2,3n})+2f(x_{1,n})\big]\,,
$$
$$
D_3(\mathbf{t})[f]=27\,n^{3}\big[-f(x_{0,n})
+3f(x_{1,3n})-3f(x_{2,3n})+f(x_{1,n})\big]\,,
$$
and $\,\widetilde{D}_1(\mathbf{t})$,
$\,\widetilde{D}_1(\mathbf{t})\,$ are obtained from
$D_1(\mathbf{t})$ and $D_3(\mathbf{t})$ by reflection. By \eqref{e9}
$D(\mathbf{t})$ generates a symmetric $(n+5)$-point quadrature
formula
$$
Q_{n+5}[f]=\sum_{k=1}^{n+5}A_{k,n+5} f(\tau_{k,n})\,,
$$
with nodes $\{\tau_{i,n+5}\}_{i=1}^{n+5}$ given by
$$
\begin{array}{l}
\tau_{k,n+5}=x_{k-1,3n},  \ 1\leq k\leq 3\vspace*{1mm}\\
\tau_{k,n+5}=x_{k-3,n}, \ \ 4\leq k\leq n+2\vspace*{1mm}\\
\tau_{k,n+5}=x_{2n-5+k,3n},\ n+3\leq k\leq n+5\,,
\end{array}
$$
and weights
\[
\begin{array}{ll}
A_{1,n+5}=A_{n+5,n+5}=\frac{43}{384n},\ \ \ &
A_{2,n+5}=A_{n+4,n+5}=\frac{69}{128n},\vspace*{2mm}\\
A_{3,n+5}=A_{n+3,n+5}=-\frac{21}{128n},\ \ \ &
A_{4,n+5}=A_{n+2,n+5}=\frac{389}{384n},\vspace*{2mm}\\
A_{k,n+5}=\frac{1}{n},\ 5\leq k\leq n+1\,. &
\end{array}
\]
\begin{figure}[h]
\begin{center}
\includegraphics[width=0.47\textwidth]{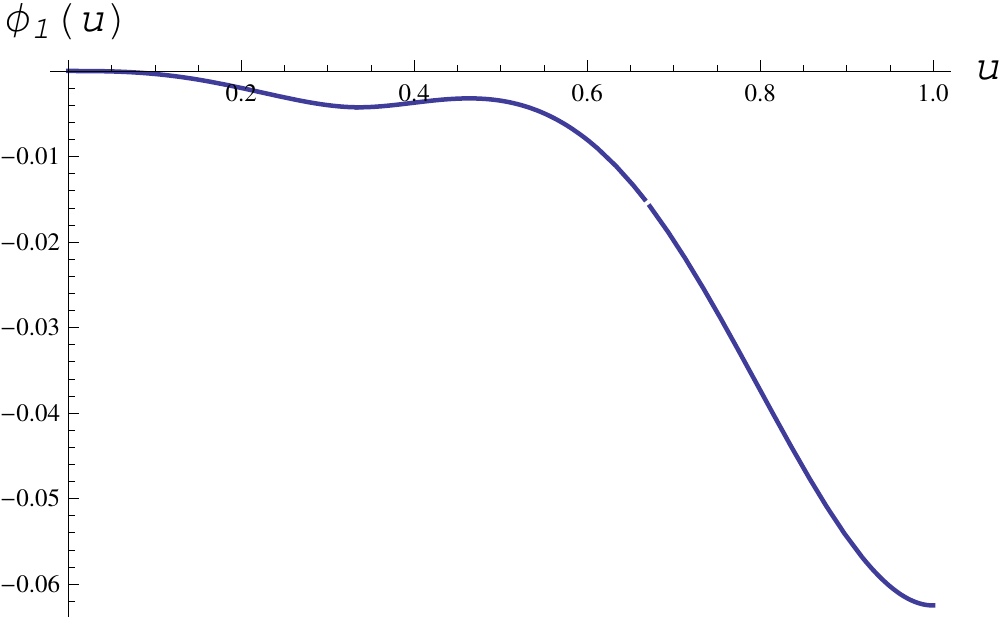}\hspace{0.05\textwidth}
\includegraphics[width=0.47\textwidth]{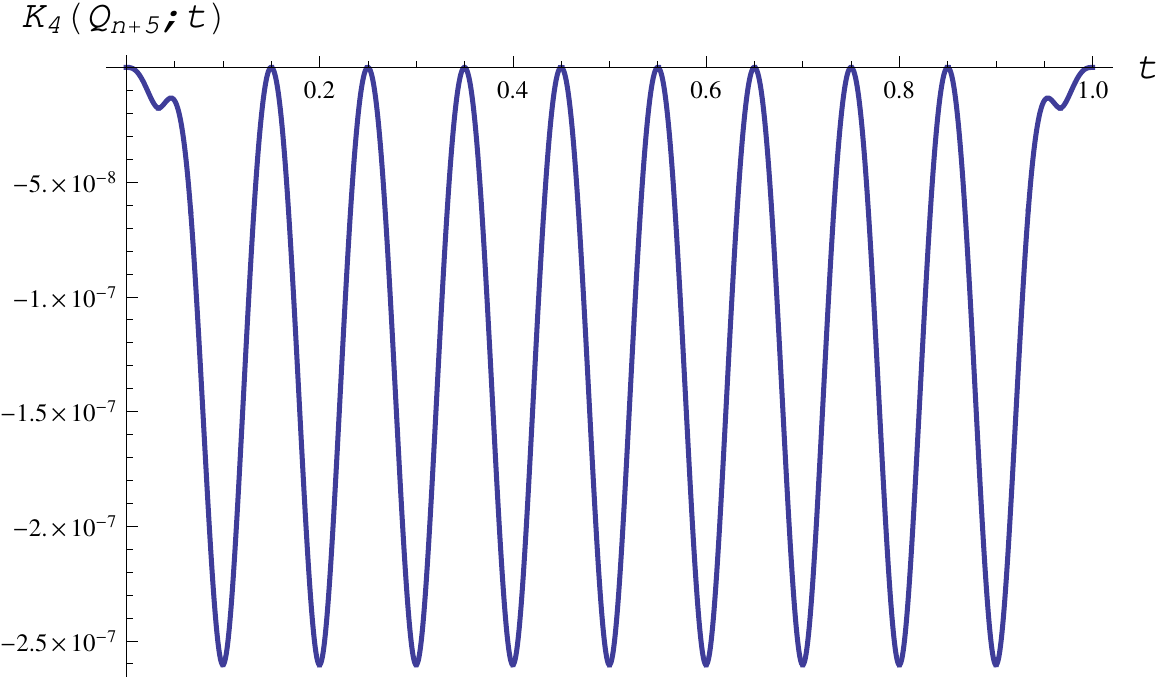}
\caption{Graphs of $\,\phi_1(u)\,$ (left), and $\,K_4(Q_{n+5};t)$,
$n=10\;$ (right).}\label{fig:1}
\end{center}
\end{figure}

In view of Proposition \eqref{p1}, to verify that $Q_{n+5}$ is a
negative definite quadrature formula of order $4$, we only have to
show that $\,K_4(Q_{n+5};t)< 0\,$ for $\,t\in (0, x_{1,n})$. By
substituting $t=u/n$, $\,u\in (0,1)$, this task reduces to
$$
\phi_1(u):=u^4-\frac{43}{96}\,u^3-\frac{69}{32}\,(u-1/3)_{+}^{3}
+\frac{21}{32}\,(u-2/3)_{+}^{3}\stackrel{?}{\leq} 0\,, \qquad u\in
(0,1)\,.
$$
The verification is straightforward, and we omit it. The graphs of
$\phi_1(u)$ and $\,K_4(Q_{n+5};t)$, with $\,n=10$, are depicted on
Figure~1.

For the evaluation of the error constant
$c_4(Q_{n+5})=I[K_4(Q_{n+5};\cdot)]$, we make use of Proposition
\ref{p1} and the fact that $K_4(Q_{n+5};\cdot)$ is symmetric, hence
\[
\begin{split}
c_4(Q_{n+5})&=2\int_{0}^{x_{1,n}}K_4(Q_{n+5};t)dt+\frac{1}{n^{4}}
\int_{x_{1,n}}^{x_{n-1,n}}\big[\widetilde{B}_4(n\,t)-B_4(1/2)\big]dt\\
&=\frac{1}{12\,n^{5}}\int_{0}^{1}\phi_1(u)\,du
-\frac{B_4(1/2)}{n^{4}}\Big(1-\frac{2}{n}\Big)\,.
\end{split}
\]
A further calculation shows that $I[\phi_1]=-71/4320$ and
$$
c(Q_{n+5})=-\frac{7}{5760\,n^4}\Big(1-\frac{55}{63\,n}\Big)\,.
$$
%=================================================================
\subsubsection{A quadrature formula generated by
$D(x_{0,n},y_{1,n},x_{1,n},x_{2,n})$}

With this set of formulae for numerical differentiation we get
through \eqref{e9} an $(n+3)$-point quadrature formula
$$
Q_{n+3}[f]=\sum_{k=1}^{n+3}A_{k,n} f(\tau_{k,n+3})
$$
with nodes
$$
\begin{array}{llll}
\tau_{1,n+3}=x_{0,n}, & \tau_{2,n+3}=y_{1,n}, &
\tau_{n+2,n+3}=y_{n-1,n}
& \tau_{n+3,n+3}=x_{n,n}\,,\vspace*{1mm}\\
\tau_{k,n+3}=x_{k-2,n}, & 3\leq k\leq n+1\,,& &
\end{array}
$$
and weights
$$
\begin{array}{ll}
A_{1,n+3}=A_{n+3,n+3}=\frac{43}{192n},\ \ \ &
A_{2,n+3}=A_{n+2,n+3}=\frac{29}{72n},\vspace*{2mm}\\
A_{3,n+3}=A_{n+1,n+3}=\frac{83}{96n},\ \ \ &
A_{4,n+3}=A_{n,n+3}=\frac{581}{576n},\vspace*{2mm}\\
A_{k,n+3}=\frac{1}{n},\ 5\leq k\leq n-1\,.&
\end{array}
$$
In view of Proposition \ref{p1}, $Q_{n+3}$ is a negative definite of
order $4$ if and only if $\,K_4(Q_{n+3};t)<0\,$ for $\,t\in
(0,x_{2,n})$. By change of variable $t=u/n$, $u\in (0,2)$, this
condition becomes
$$
\phi_2(u):=u^{4}-\frac{43}{48}\,u^{3}-\frac{29}{18}\,(u-1/2)_{+}^3
-\frac{83}{24}\,(u-1)_{+}^3<0\,,\quad u\in (0,2)\,,
$$
and it is not difficult to verify that it is fulfilled.
%\begin{figure}
%\begin{center}
%\includegraphics[width=0.47\textwidth]{fig-2a.pdf}\hspace{0.05\textwidth}
%\includegraphics[width=0.47\textwidth]{fig-2b.pdf}
%\caption{Graphs of $\,\phi_2(u)\,$ (left), and $\,K_4(Q_{n+3};t)$,
%$\,n=10\;$ (right).}\label{fig:2}
%\end{center}
%\end{figure}

For the error constant $c_4(Q_{n+3})=I[K_4(Q_{n+3};\cdot)]$ we have
\[
\begin{split}
c_4(Q_{n+3})&=2\int_{0}^{x_{2,n}}K_4(Q_{n+3};t)dt+\frac{1}{n^{4}}
\int_{x_{2,n}}^{x_{n-2,n}}\big[\widetilde{B}_4(n\,t)-B_4(1/2)\big]dt\\
&=\frac{1}{12\,n^{5}}\int_{0}^{2}\phi_2(u)\,du
-\frac{B_4(1/2)}{n^{4}}\Big(1-\frac{4}{n}\Big)\,,
\end{split}
\]
and after evaluation of the integral of $\phi_2$, we obtain
$$
c(Q_{n+3})=-\frac{7}{5760\,n^4}\Big(1+\frac{55}{28\,n}\Big)\,.
$$

%================================================================
\subsection{Negative definite quadrature formulae of order $4$ based on $Q_n^{Mi}$}
We rewrite the first Euler-Maclaurin summation formula in Lemma
\ref{l1} with $s=4$ in the form
\[
\begin{split}
\int\limits_{0}^{1}\!\!f(x)dx&\!=\! Q_{n}^{Mi}[f]\!
+\!\frac{1}{24\,n^2}\big[f^{\prime}(1)\!-\!f^{\prime}(0)\big]\!+\!
\frac{1}{n^{4}}\!\int\limits_0^{1}\!\!
\big[\widetilde{B}_4(nx-\frac{1}{2})-B_4(\frac{1}{2})\big]\,f^{(4)}(x)dx\\
&\!=:Q_n^{\prime\prime}[f]+\frac{1}{n^{4}}\!\int\limits_0^{1}\!\!
\big[\widetilde{B}_4(nx-\frac{1}{2})-B_4(\frac{1}{2})\big]\,f^{(4)}(x)dx\,.
\end{split}
\]
Here, $Q_n^{\prime\prime}$ is a negative definite quadrature formula
of order $4$, as, by \eqref{e7}, its fourth Peano kernel
$K_4(Q_n^{\prime\prime};x)=n^{-4}\big[\widetilde{B}_4(nx-1/2)-B_4(1/2)\big]\,$
is non-positive. Since $Q_n^{\prime\prime}$ is not of the desired
form, we choose a set $D(\mathbf{t})$ of formulae for numerical
differentiation to replace the values of $f^{\prime}$ in
$Q_n^{\prime\prime}[f]$, and thus to obtain a (symmetric) quadrature
formula
\begin{equation}\label{e10}
Q=Q_n^{Mi}-\frac{1}{12\,n^2}\big(D_1(\mathbf{t})-
\widetilde{D}_1(\mathbf{t})\big)\,,
\end{equation}
which involves at most $8$ nodes in addition to
$\,\{y_{k,n}\}_{k=1}^{n}$. By the same argument that led us to
Proposition \ref{p1}, here we have
\begin{prop}\label{p2}
The fourth Peano kernel of the symmetric quadrature formula $Q$
generated by $D(\mathbf{t})$ through \eqref{e10} satisfies
$$
K_4(Q;x)\equiv \frac{1}{n^{4}}\,
\big[\widetilde{B}_4(nx-\frac{1}{2})-B_4(\frac{1}{2})\big]\,, \qquad
x\in [t_4,1-t_4]\,.
$$
Consequently, $\,Q\,$ is negative definite of order $4$ if and only
if $\,K_4(Q;x)\leq 0\,$ for $\,x\in (0,t_4)$\,.
\end{prop}
%===================================================================
\subsubsection{A quadrature formula generated by
$D(x_{0,n},y_{1,n},y_{2,2n},x_{1,n})$} For
$\,\mathbf{t}=(x_{0,n},y_{1,n},y_{2,2n},x_{1,n})\,$ we have
$$
D_{1}(\mathbf{t})[f]=\frac{n}{3}
\bigg[-13f(x_{0,n})+36f(y_{1,n})-32f(x_{3,4n})+9f(x_{1,n})\bigg],
$$
and $D(\mathbf{t})$ generates through \eqref{e10} a $(n+6)$-point
quadrature formula
$$
Q_{n+6}[f]=\sum_{k=1}^{n+6}A_{k,n+6}\,f(\tau_{k,n+6})
$$
with nodes
$$
\begin{array}{l}
\tau_{1,n+6}=x_{0,n},\ \ \tau_{2,n+6}=y_{1,n},\ \
\tau_{3,n+6}=y_{2,2n},\ \ \tau_{4,n+6}=x_{1,n},\vspace*{1mm}\\
\tau_{k,n+6}=y_{k-3,n},\ \ \ \ \ 5\leq k\leq n+2,\vspace*{1mm}\\
\tau_{n+7-k,n+6}=1-\tau_{k,n+6},\ \ \ \ 1\leq k\leq 4
\end{array}
$$
and weights
$$
\begin{array}{ll}
A_{1,n+6}=A_{n+6,n+6}=\frac{13}{72n},\ \ &
A_{2,n+6}=A_{n+5,n+6}=\frac{1}{2n},\vspace*{2mm}\\
A_{3,n+6}=A_{n+4,n+6}=\frac{4}{9n},\ &
A_{4,n+6}=A_{n+3,n+6}=-\frac{1}{8n},\vspace*{2mm}\\
A_{k,n+6}=\frac{1}{n},\ \ \ 5\leq k\leq n+2\,.&
\end{array}
$$

By Proposition \ref{p2}, to verify that $Q_{n+6}$ is a negative
definite quadrature formula of order $4$, we only need to check
whether $\,K_4(Q_{n+6};t)<0,\,$ $\,t\in (0,x_{1,n})$, which, after
the change of variable $t=u/n$, becomes
$$
\phi_3(u):=u^{4}-\frac{13}{18}\,u^{3}-2\,(u-1/2)_{+}^3
-\frac{16}{9}\,(u-3/4)_{+}^3\stackrel{?}{\leq} 0\,,\quad u\in
(0,1)\,.
$$
\begin{figure}[h]
\begin{center}
\includegraphics[width=0.47\textwidth]{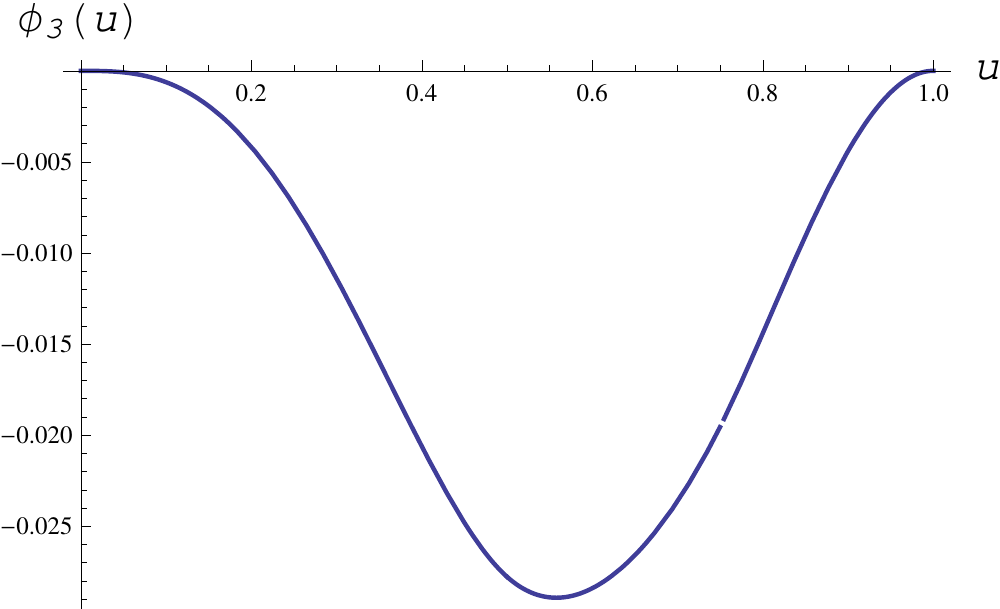}\hspace{0.05\textwidth}
\includegraphics[width=0.47\textwidth]{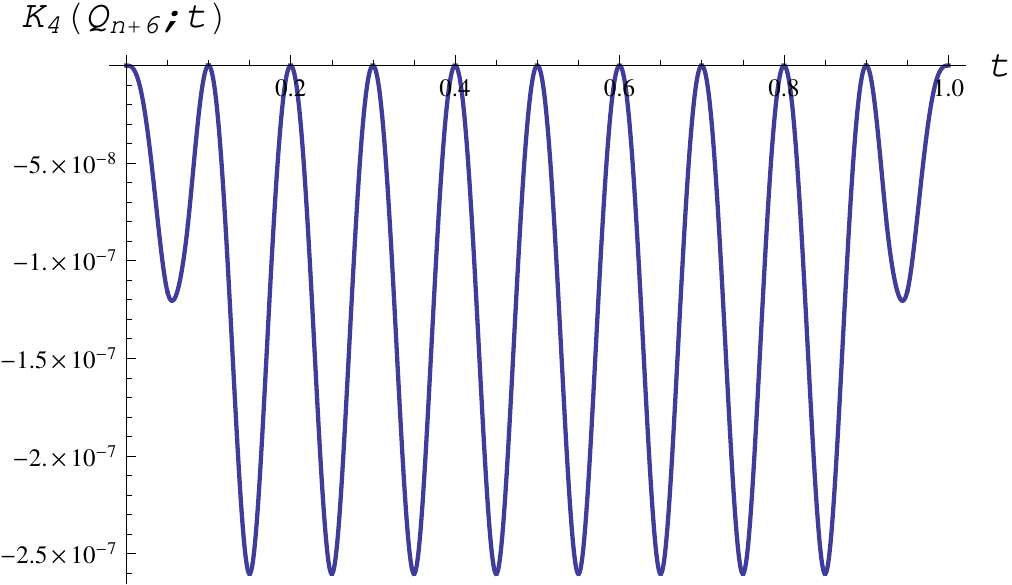}
\caption{Graphs of $\,\phi_3(u)\,$ (left), and $\,K_4(Q_{n+6};t)$,
$\,n=10\;$ (right).}\label{fig:2}
\end{center}
\end{figure}

The latter condition is fulfilled, as is seen also on Fugure
\ref{fig:2} (left). For the error constant $c_4(Q_{n+6})$, in view
of Proposition \ref{p2}, we have
\[
\begin{split}
c_4(Q_{n+6})&=2\int_{0}^{x_{1,n}}K_4(Q_{n+6};t)dt+\frac{1}{n^{4}}
\int_{x_{1,n}}^{x_{n-1,n}}\big[\widetilde{B}_4(n\,t-1/2)-B_4(1/2)\big]dt\\
&=\frac{1}{12\,n^{5}}\int_{0}^{1}\phi_3(u)\,du
-\frac{B_4(1/2)}{n^{4}}\Big(1-\frac{2}{n}\Big)\,.
\end{split}
\]
With further calculations we find $I[\phi_3]=-13/960$ and
$$
c_4(Q_{n+6})=-\frac{7}{5760\,n^4}\Big(1-\frac{15}{14\,n}\Big)\,.
$$

Below we give two further negative definite quadrature formulae
generated through \eqref{e10}.
%===================================================================
\subsubsection{A quadrature formula generated by
$D(x_{0,n},y_{1,2n},y_{1,n},x_{1,n})$} For
$t=(x_{0,n},y_{1,2n},y_{1,n},x_{1,n})$ we have
$$
D_{1}(\mathbf{t})[f]=\frac{n}{3}
\bigg[-21f(x_{0,n})+32f(y_{1,2n})-12f(y_{1,n})+f(x_{1,n})\bigg],
$$
and $D(\mathbf{t})$ generates through \eqref{e10} another
$(n+6)$-point quadrature formula
$$
Q_{n+6}[f]=\sum_{k=1}^{n+6}A_{k,n+6}\,f(\tau_{k,n+6})
$$
with nodes
$$
\begin{array}{l}
\tau_{1,n+6}=x_{0,n},\ \ \tau_{2,n+6}=y_{1,2n},\ \
\tau_{3,n+6}=y_{1,n},\ \ \tau_{4,n+6}=x_{1,n},\vspace*{1mm}\\
\tau_{k,n+6}=y_{k-3,n},\ \ \ \ \ \ \ 5\leq k\leq n+2\,,\vspace*{1mm}\\
\tau_{n+7-k,n+6}=1-\tau_{k,n+6}\,,\vspace*{1mm}\ \ \ \ \ \ 1\leq
k\leq 4\,
\end{array}
$$
and weights
$$
\begin{array}{ll}
A_{1,n+6}=A_{n+6,n+6}=\frac{7}{24n},\ \ \ &
A_{2,n+6}=A_{n+5,n+6}=-\frac{4}{9n},\vspace*{2mm}\\
A_{3,n+6}=A_{n+4,n+6}=\frac{7}{6n},\ \ \ &
A_{4,n+6}=A_{n+3,n+6}=-\frac{1}{72n},\vspace*{2mm}\\
A_{k,n+6}=\frac{1}{n},\ \ \ 5\leq k\leq n+2\,.&
\end{array}
$$
The error constant of $Q_{n+6}$ is
$$
c_4(Q_{n+6})=-\frac{7}{5760\,n^4}\Big(1-\frac{5}{14\,n}\Big)\,.
$$
%==================================================================
\subsubsection{A quadrature formula generated by $D(x_{0,n},
y_{1,6n}, y_{1,3n}, y_{1,2n})$} With $\mathbf{t}=(x_{0,n}, y_{1,6n},
y_{1,3n}, y_{1,2n})$, $\,D(\mathbf{t})$ generates a negative
definite of order~$4\,$ $\,(n+8)$-point quadrature formula
$$
Q_{n+8}[f]=\sum_{k=1}^{n+8}A_{k,n+8}\,f(\tau_{k,n+8})
$$
with nodes, weights and error constant given by
$$
\begin{array}{l}
\tau_{1,n+8}=x_{0,n},\ \ \tau_{2,n+8}=y_{1,6n},\ \
\tau_{3,n+8}=y_{1,3n},\ \ \tau_{4,n+8}=y_{1,2n},\vspace*{1mm}\\
\tau_{k,n+8}=y_{k-4,n}, \ \ \ \ \ \ \ \ 5\leq k\leq n+4,\vspace*{1mm}\\
\tau_{n+9-k,n+8}=1-\tau_{k,n+8}, \ \ \ \ \ \ \ 1\leq k\leq 4,
\end{array}
$$
$$
\begin{array}{ll}
A_{1,n+8}=A_{n+8,n+8}=\frac{11}{12n},\ \ \ &
A_{2,n+8}=A_{n+7,n+8}=-\frac{3}{2n},\vspace*{2mm}\\
A_{3,n+8}=A_{n+6,n+8}=\frac{3}{4n},\ \ \ &
A_{4,n+8}=A_{n+5,n+8}=-\frac{1}{6n},\vspace*{2mm}\\
A_{k,n+8}=\frac{1}{n},\ \ \ \ \ \ 5\leq k\leq n+4\,,&
\end{array}
$$
$$
c_4(Q_{n+8})=-\frac{7}{5760\,n^4}\Big(1-\frac{5}{504\,n}\Big)\,.
$$
Clearly, the error constant $c_4(Q_{n+8})$ is inferior to those of
the preceding two quadrature formulae, which moreover involve two
nodes less. The reason for quoting this quadrature formula will
become clear in Section 4.
%==================================================================
\subsection{Positive definite quadrature formulae of order $4$ based on
$Q_n^{Tr}$} We rewrite the second Euler-Maclaurin summation formula
in Lemma \ref{l1} with $s=4$ in the form
\[
\begin{split}
\int\limits_{0}^{1}f(x)\,dx=&Q_{n}^{Tr}[f]
-\frac{1}{12\,n^2}\big[f^{\prime}(1)-f^{\prime}(0)\big]
+\frac{1}{n^{4}}\int\limits_0^{1}
\big[\widetilde{B}_4(nx)-B_4(0)\big]\,f^{(4)}(x)\,dx\\
=:&\widetilde{Q}_n^{\prime}[f]+\frac{1}{n^{4}}\int\limits_0^{1}
\big[\widetilde{B}_4(nx)-B_4(0)\big]\,f^{(4)}(x)\,dx\,.
\end{split}
\]
By \eqref{e7}, $\widetilde{Q}_n^{\prime}$ is a positive definite
quadrature formula of order $4$, and we choose a set of formulae for
numerical differentiation $D(\mathbf{t})$ to approximate
$f^{\prime}(0)$ and $f^{\prime}(1)$ in $\widetilde{Q}_n^{\prime}$,
thus arriving at a new quadrature formula
\begin{equation}\label{e11}
Q=Q_n^{Tr}+\frac{1}{12\,n^2}\,
\big(D_1(\mathbf{t})-\widetilde{D}_1(\mathbf{t})\big)\,,
\end{equation}
which involves at most $8$ nodes in addition to
$\{x_{k,n}\}_{k=0}^{n}$\,.

A statement analogous to Propositions \ref{p1}, \ref{p2} holds true:
\begin{prop}\label{p3}
The fourth Peano kernel of the symmetric quadrature formula $Q$
generated by $D(\mathbf{t})$ through \eqref{e11} satisfies
$$
K_4(Q;x)\equiv \frac{1}{n^{4}}\,
\big[\widetilde{B}_4(nx)-B_4(0)\big]\,, \qquad x\in [t_4,1-t_4]\,.
$$
Consequently, $\,Q\,$ is positive definite of order $4$ if and only
if $\,K_4(Q;x)\geq 0\,$ for $\,x\in (0,t_4)$\,.
\end{prop}

Below we construct three positive definite quadrature formulae
generated through \eqref{e11} by different sets $D(\mathbf{t})$ of
formulae for numerical differentiation.
%========================================================
\subsubsection{A quadrature formula generated by
$D(x_{0,n},y_{1,3n},x_{1,3n},y_{1,n})$} With
$\mathbf{t}=(x_{0,n},y_{1,3n},x_{1,3n},y_{1,n})$, $D(\mathbf{t})$
generates through \eqref{e11} an $(n+7)$-point symmetric quadrature
formula
$$
Q_{n+7}[f]=\sum_{k=1}^{n+7}A_{k,n+7}\,f(\tau_{k,n+7})
$$
with nodes and weights given by
$$
\begin{array}{l}
\tau_{1,n+7}=x_{0,n},\ \ \tau_{2,n+7}=y_{1,3n},\ \
\tau_{3,n+7}=x_{1,3n},\ \ \tau_{4,n+7}=y_{1,n},\vspace*{1mm}\\
\tau_{k,n+7}=x_{k-4,n}\,,\ \ \ \ \ \ 5\leq k\leq n+3\,,\vspace*{1mm}\\
\tau_{n+8-k,n+7}=1-\tau_{k,n+7}\,, \ \ \ \ \ \ 1\leq k\leq 4\,,
\end{array}
$$
$$
\begin{array}{ll}
A_{1,n+7}=A_{n+7,n+7}=-\frac{5}{12n},\ \ \ &
A_{2,n+7}=A_{n+6,n+7}=\frac{3}{2n},\vspace*{2mm}\\
A_{3,n+7}=A_{n+5,n+7}=-\frac{3}{4n},\ \ \ &
A_{4,n+7}=A_{n+4,n+7}=\frac{1}{6n},\vspace*{2mm}\\
A_{k,n+7}=\frac{1}{n},\ \ \ 5\leq k\leq n+3\,.&
\end{array}
$$
By Proposition \ref{p3}, the verification that $\,Q_{n+7}\,$ is
positive definite of order $4$ reduces to $\,K_4(Q_{n+7};t)\geq 0\,$
for $\,t\in (0,y_{1,n})$, which, after the change of variable
$t=u/n$, $u\in (0,1/2)$, becomes
$$
\psi(u):=u^4+\frac{5}{3}\,u^3-6(u-1/6)_{+}^{3}+3(u-1/6)_{+}^{3}
\stackrel{?}{\geq} 0\,,\quad u\in (0,1/2)\,.
$$
\begin{figure}[h]
\begin{center}
\includegraphics[width=0.47\textwidth]{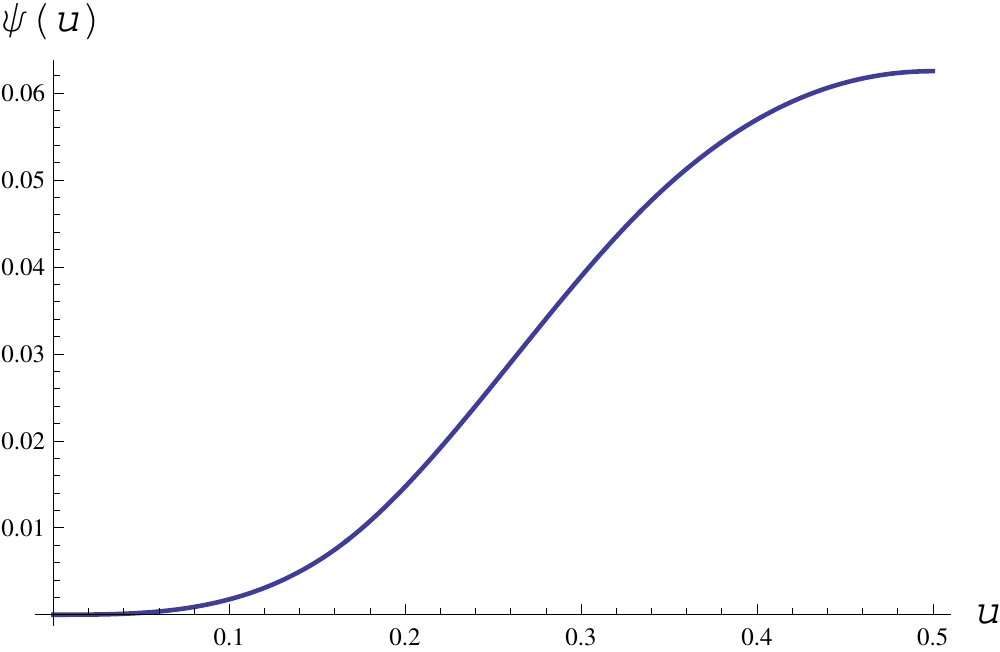}\hspace{0.05\textwidth}
\includegraphics[width=0.47\textwidth]{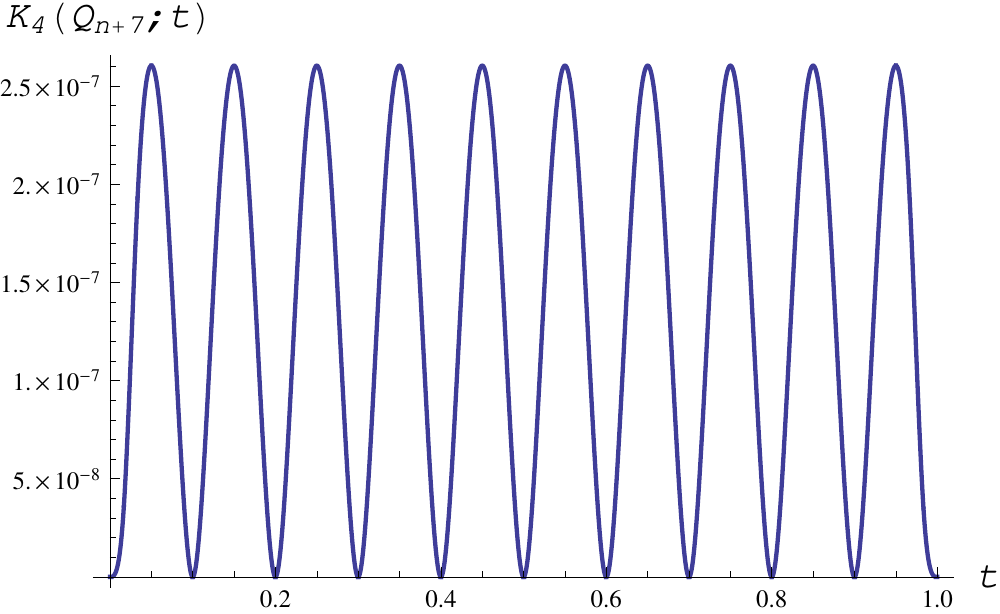}
\caption{Graphs of $\,\psi(u)\,$ (left), and $\,K_4(Q_{n+7};t)$,
$\,n=10\;$ (right).}\label{fig:3}
\end{center}
\end{figure}

The graph of $\psi$, depicted on Fugure \ref{fig:3} (left), shows
that, indeed $\psi(u)>0$ for $u\in (0,1/2)$. Finally, in view of
Proposition \ref{p3}, for the error constant $c_4(Q_{n+7})$ we have
\[
\begin{split}
c_4(Q_{n+7})&=2\int_{0}^{y_{1,n}}K_4(Q_{n+6};t)dt+\frac{1}{n^{4}}
\int_{y_{1,n}}^{y_{n-1,n}}\big[\widetilde{B}_4(n\,t)-B_4(0)\big]dt\\
&=\frac{1}{12\,n^{5}}\int_{0}^{1/2}\psi(u)\,du
-\frac{B_4(0)}{n^{4}}\Big(1-\frac{1}{n}\Big)\,.
\end{split}
\]
The integral of $\psi$ is equal to $31/2160$, and a further
simplification implies that
$$
c_4(Q_{n+7})=\frac{1}{720\,n^{4}}\,\Big(1-\frac{5}{36\,n}\Big).
$$

The next two quadrature formulae are obtained through the same
scheme. We only give their nodes, weights and error constants,
skipping the details on the verification of their definiteness and
the calculations, as these go along the same lines as in the case we
just considered.
%================================================================
\subsubsection{A quadrature formula generated by
$D(x_{0,n},y_{1,2n},y_{1,n},x_{1,n})$} With
$\mathbf{t}=(x_{0,n},y_{1,2n},y_{1,n},x_{1,n})$, $D(\mathbf{t})$
generates through \eqref{e11} an $(n+5)$-point symmetric quadrature
formula
$$
Q_{n+5}[f]=\sum_{k=1}^{n+5}A_{k,n+5}\,f(\tau_{k,n+5})
$$
with nodes and weights given by
$$
\begin{array}{l}
\tau_{1,n+5}=x_{0,n},\ \ \tau_{2,n+5}=y_{1,2n},\ \
\tau_{3,n+5}=y_{1,n},\ \ \tau_{4,n+5}=x_{1,n}, \vspace*{1mm}\\
\tau_{k,n+5}=x_{k-3,n},\ \ \ \ \ \ \  5\leq k\leq n+1, \vspace*{1mm}\\
\tau_{n+6-k,n+5}=1-\tau_{k,n+5},\ \ \ \ 1\leq k\leq 4\,,
\end{array}
$$
$$
\begin{array}{ll}
A_{1,n+5}=A_{n+5,n+5}=-\frac{1}{12n},\ \ \ &
A_{2,n+5}=A_{n+4,n+5}=\frac{8}{9n},\vspace*{2mm}\\
A_{3,n+5}=A_{n+3,n+5}=-\frac{1}{3n},\ \ \ &
A_{4,n+5}=A_{n+2,n+5}=\frac{37}{36n},\vspace*{2mm}\\
A_{k,n+5}=\frac{1}{n},\ \ \ \ \ \ \ 5\leq k\leq n+1\,.&
\end{array}
$$
The error constant of $Q_{n+5}$ is
$$
c_4(Q_{n+5})=\frac{1}{720\,n^4}\,\Big(1-\frac{5}{8\,n}\Big)\,.
$$
%=============================================================
\subsubsection{A quadrature formula generated by
$D(x_{0,n},y_{1,2n},y_{1,n},y_{2,2n})$} With
$\mathbf{t}=(x_{0,n},y_{1,2n},y_{1,n},y_{2,n})$, $D(\mathbf{t})$
generates through \eqref{e11} an $(n+7)$-point symmetric quadrature
formula
$$
Q_{n+7}[f]=\sum_{k=1}^{n+7}A_{k,n+7}\,f(\tau_{k,n+7})
$$
with nodes, weights and error constant given by
$$
\begin{array}{l}
\tau_{1,n+7}=x_{0,n},\ \ \tau_{2,n+7}=y_{1,2n},\ \
\tau_{3,n+7}=y_{1,n},\ \ \tau_{4,n+7}=y_{2,2n},\vspace*{1mm}\\
\tau_{k,n+7}=x_{k-4,n}\,,\ \ \ \ \ \ 5\leq k\leq n+3\,,\vspace*{1mm}\\
\tau_{n+8-k,n+7}=1-\tau_{k,n+7}\,, \ \ \ \ \ \ 1\leq k\leq 4\,,
\end{array}
$$
$$
\begin{array}{ll}
A_{1,n+7}=A_{n+7,n+7}=-\frac{1}{9n},\ \ \ &
A_{2,n+7}=A_{n+6,n+7}=\frac{1}{n},\vspace*{2mm}\\
A_{3,n+7}=A_{n+5,n+7}=-\frac{1}{2n},\ \ \ &
A_{4,n+7}=A_{n+4,n+7}=\frac{1}{9n},\vspace*{2mm}\\
A_{k,n+7}=\frac{1}{n},\ \ \ 5\leq k\leq n+3\,,&
\end{array}
$$
$$
c_4(Q_{n+7})=\frac{1}{720\,n^{4}}\,\Big(1-\frac{15}{32\,n}\Big).
$$
%==================================================================
\subsection{Positive definite quadrature formulae of order $4$ based on
$Q_n^{Mi}$} We write the first formula in Lemma \ref{l1} with $s=4$
in the form
\[
\begin{split}
\int\limits_{0}^{1}f(x)\,dx=&Q_{n}^{Mi}[f]+\frac{1}{24\,n^2}\big[f^{\prime}(1)-f^{\prime}(0)\big]-
\frac{1}{384\,n^4}\big[f^{\prime\prime\prime}(1)-f^{\prime\prime\prime}(0)\big]\\
&+\frac{1}{n^{4}}\int\limits_{0}^{1}
\Big[\widetilde{B}_4\Big(nx-\frac{1}{2}\Big)-B_4(0)\Big]\,f^{(4)}(x)\,dx\\
=:&\widetilde{Q}_n^{\prime\prime}[f]+
\frac{1}{n^{4}}\int\limits_{0}^{1}
\Big[\widetilde{B}_4\Big(nx-\frac{1}{2}\Big)-B_4(0)\Big]\,f^{(4)}(x)\,dx\,.
\end{split}
\]
We choose a set $D(\mathbf{t})$ for approximating the derivatives
values in $\widetilde{Q}_n^{\prime\prime}$, thus obtaining a
quadrature formula
\begin{equation}\label{e12}
Q=Q_n^{Mi}-\frac{1}{24\,n^2}\big(D_1(\mathbf{t})-
\widetilde{D}_1(\mathbf{t})\big)+\frac{1}{384\,n^{4}}\big(D_3(\mathbf{t})
-\widetilde{D}_3(\mathbf{t})\big)\,,
\end{equation}
which involves at most $8$ nodes in addition to
$\{y_{i,n}\}_{i=1}^{n}$. We have
\begin{prop}\label{p4}
The fourth Peano kernel of the symmetric quadrature formula $Q$
generated by $D(\mathbf{t})$ through \eqref{e12} satisfies
$$
K_4(Q;x)\equiv \frac{1}{n^{4}}\,
\big[\widetilde{B}_4(nx-\frac{1}{2})-B_4(0)\big]\,, \qquad x\in
[t_4,1-t_4]\,.
$$
Consequently, $\,Q\,$ is positive definite of order $4$ if and only
if $\,K_4(Q;x)\geq 0\,$ for $\,x\in (0,t_4)$\,.
\end{prop}

On using Proposition \ref{p4}, we verify the definiteness  and
evaluate the error constant of $Q$. We give below two positive
definite quadrature formulae of order $4$, constructed on the basis
of \eqref{e12}. As the definiteness verification and the evaluation
of the error constants are completely analogous to that in the
preceding cases, they are skipped here.
%==================================================================
\subsubsection{A quadrature formula generated by
$D(y_{1,n},x_{1,n},y_{2,n},y_{3,n})$} With
$\mathbf{t}=(y_{1,n},x_{1,n},y_{2,n},y_{3,n})$, $\,n\geq 7$, we
obtain through \eqref{e12} an $(n+2)$-point symmetric quadrature
formula
$$
Q_{n+2}[f]=\sum_{k=1}^{n+2}A_{k,n+2}\,f(\tau_{k,n+2})\,,
$$
which is positive definite of order $4$. The nodes and the weights
 of $Q_{n+2}$ are
$$
\begin{array}{l}
\tau_{1,n+2}=y_{1,n},\ \ \tau_{2,n+2}=x_{1,n},\ \
\tau_{n+1,n+2}=x_{n-1,n},\ \  \tau_{n+2,n+2}=y_{n,n}
\,,\vspace*{1mm}\\
\tau_{k,n+2}=y_{k-1,n},\ \ \ \ \ \ \ 3\leq k\leq n\,.
\end{array}
$$
$$
\begin{array}{ll}
A_{1,n+2}=A_{n+2,n+2}=\frac{251}{192n},\ \ \ &
A_{2,n+2}=A_{n+1,n+2}=-\frac{43}{72n},\vspace*{2mm}\\
A_{3,n+2}=A_{n,n+2}=\frac{127}{96n},\ \ \ &
A_{4,n+2}=A_{n-1,n+2}=\frac{557}{576n},\vspace*{2mm}\\
A_{k,n+2}=\frac{1}{n},\ \ \ \ \ \ \ 5\leq k\leq n-2\,.&
\end{array}
$$
The error constant of $Q_{n+2}$ is
$$
c_4(Q_{n+2})=\frac{1}{720\,n^4}\,\Big(1+\frac{445}{32\,n}\Big)\,.
$$
%==================================================================
\subsubsection{A quadrature formula generated by
$D(x_{0,n},y_{1,3n},x_{1,3n},y_{1,n})$} With
$\mathbf{t}=(x_{0,n},y_{1,3n},x_{1,3n},y_{1,n})$, $\,n\geq 3$, we
obtain through \eqref{e12} an $(n+6)$-point positive definite of
order $4$ quadrature formula
$$
Q_{n+6}[f]=\sum_{k=1}^{n+6}A_{k,n+6}\,f(\tau_{k,n+6})\,,
$$
with nodes, weights and error constant given by
$$
\begin{array}{l}
\tau_{1,n+6}=x_{0,n},\ \ \tau_{2,n+6}=y_{1,3n},\ \
\tau_{3,n+6}=x_{1,3n},\ \ \tau_{4,n+6}=y_{1,n},\vspace*{1mm}\\
\tau_{k,n+6}=y_{k-3,n},\ \ \ \ \ \ \ \ \ \ 5\leq k\leq n+2,\vspace*{1mm}\\
\tau_{n+7-k,n+6}=\tau_{k,n+6},\ \ \ \ 1\leq k\leq 4\,,
\end{array}
$$
$$
\begin{array}{ll}
A_{1,n+6}=A_{n+6,n+6}=-\frac{5}{48n},\ \ \ &
A_{2,n+6}=A_{n+5,n+6}=\frac{15}{16n},\vspace*{2mm}\\
A_{3,n+6}=A_{n+4,n+6}=-\frac{21}{16n},\ \ \ &
A_{4,n+6}=A_{n+3,n+6}=\frac{71}{48n},\vspace*{2mm}\\
A_{k,n+6}=\frac{1}{n},\ \ \ \ \ \ \ 5\leq k\leq n+2\,,&
\end{array}
$$
$$
c_4(Q_{n+6})=\frac{1}{720\,n^4}\,\Big(1-\frac{125}{144\,n}\Big)\,.
$$
%===================================================================
\section{Monotonicity of the remainders and a posteriori error
estimates}

In this section we shall exploit the following general observation
about definite quadrature formulae.
\begin{thm}\label{t1}
Let $\,(Q^{\prime}\,,\,Q^{\prime\prime})\,$ be a pair of positive
(negative) definite quadrature formulae of order $r$. Assume that,
for some $c>0$, the quadrature formula
$$
\widehat{Q}:=(c+1)\,Q^{\prime}-c\,Q^{\prime\prime}
$$
is negative (positive) definite of order $r$. Then the following
inequalities hold true whenever $\,f\,$ is an $\,r$-convex or
$r$-concave function:
\begin{enumerate}[~~(i)~~]
\item
$\ds{|R[Q^{\prime};f]|\leq\frac{c}{c+1}\,|R[Q^{\prime\prime};f]|}$\,;
\item $\ds{|R[Q^{\prime};f]|\leq
c\,|Q^{\prime}[f]-Q^{\prime\prime}[f]|}$\,;
\item $\ds{|R[Q^{\prime\prime};f]|\leq
(c+1)\,|Q^{\prime}[f]-Q^{\prime\prime}[f]|}$\,.
\end{enumerate}
\end{thm}
\begin{pf}
Let us consider, e.g., the case when $Q^{\prime}$ and
$Q^{\prime\prime}$ are negative definite and $\,\widehat{Q}\,$ is
positive definite, of order $r$. Without loss of generality we may
assume that $f$ is $r$-convex. Then $\,R[Q^{\prime};f]\leq 0$,
$\,R[Q^{\prime\prime};f]\leq 0$, and $\,R[\widehat{Q};f]\geq 0$,
therefore
$$
0\leq
R[\widehat{Q};f]=(c+1)\,R[Q^{\prime};f]-c\,R[Q^{\prime\prime};f]\,,
$$
and hence
$$
-R[Q^{\prime};f]\leq-\frac{c}{c+1}\,(R[Q^{\prime\prime};f])\,,
$$
which, in this case, is exactly claim (i) of Theorem \ref{t1}. Claim
(iii) follows from
\[
\begin{split}
|Q^{\prime}[f]-Q^{\prime\prime}[f]|&=|R[Q^{\prime\prime};f]-R[Q^{\prime};f]|
\geq |R[Q^{\prime\prime};f]|-|R[Q^{\prime};f]|\\
&\geq |R[Q^{\prime\prime};f]|-\frac{c}{c+1}\,|R[Q^{\prime\prime};f]|
=\frac{1}{c+1}\,|R[Q^{\prime\prime};f]|\,,
\end{split}
\]
and (ii) is a consequence of (iii) and (i). The proof of the case
when $Q^{\prime}$ and $Q^{\prime\prime}$ are positive definite and
$\,\widehat{Q}\,$ is negative definite of order $r$ is analogous,
and we omit it.\qed
\end{pf}

\begin{rmk}\label{r1}
Notice the non-symmetric roles of $Q^{\prime}$ and
$Q^{\prime\prime}$ in Theorem \ref{t1}. Part~(i) implies that for
$r$-convex (concave) integrand $f$, $\,Q^{\prime}[f]\,$ furnishes a
better approximation to $I[f]$ than $\,Q^{\prime\prime}[f]$. Another
observation is that, the smaller $c>0$, the better a posteriori
error estimates (ii) and (iii) we get. Hence, it makes sense to
search for the best possible (i.e., the smallest) $\,c>0\,$ for
which $\widehat{Q}$ is definite with the opposite type of
definiteness to those of $Q^{\prime}$ and $Q^{\prime\prime}$.
\end{rmk}

\begin{ex}
If $\,(Q^{\prime},Q_n^{\prime\prime})=(Q_{2n}^{Tr},Q_n^{Tr})$, then,
since $\,\widehat{Q}=2Q_{2n}^{Tr}-Q_n^{Tr}=Q_{2n}^{Mi}\,$, the
assumptions of Theorem \ref{t1} are fulfilled with $r=2$ and $c=1$.
Hence, for $f$ convex, we have the (well-known) inequalities: $
\ds{|R[Q_{2n}^{Tr}[f]|\leq \frac{1}{2}\,|R[Q_{n}^{Tr}[f]|}$,
$\,\ds{|R[Q_{2n}^{Tr}[f]|\leq |Q_{n}^{Tr}[f]-Q_{2n}^{Tr}[f]|}$, and
$\ds{|R[Q_{n}^{Tr}[f]|\leq 2\,|Q_{n}^{Tr}[f]-Q_{2n}^{Tr}[f]|}$\,.
\end{ex}

Theorem \ref{t1} is applicable to some pairs of the definite
quadrature formulae of order $4$, obtained in Section 3. In Tables 1
and 2 below, the notation $Q_{(3.b.c),m}$ stands for the quadrature,
given in Section 3.b.c, with a parameter $n=m$.

\begin{thm}\label{t2} The assumptions of Theorem \ref{t1} are
fulfilled for the pairs $\,(Q^{\prime}, Q^{\prime\prime})\,$ of
negative definite quadrature formulae and with the best possible
constants $c$, given in Table~1.
\end{thm}

\begin{table}[h]\label{tab:1}
\begin{center}
\begin{tabular}{|c|c|c|c|}
\hline No. & $\phantom{\begin{array}{@{}l@{}}1\\[-1ex]2 \end{array}}$\ \ \ \  $Q^{\prime}$\ \ \ \ \ &\ \ \ \ \
$Q^{\prime\prime}$\ \ \ \ &\ \ \ \ $c$\ \ \ \ \\
\hline $1$ &
$\phantom{\begin{array}{@{}l@{}}1\\[-1ex]2 \end{array}}$ $Q_{(3.2.1),2n}$\ \ &\  $Q_{(3.1.1),n}$ & $\frac{104}{299}$\\
$2$ &
$\phantom{\begin{array}{@{}l@{}}1\\[-1.6ex]2 \end{array}}$ $Q_{(3.2.1),2n}$\ \ &\  $Q_{(3.1.3),n}$ & $\frac{52}{77}$\\
$3$ &
$\phantom{\begin{array}{@{}l@{}}1\\[-1.6ex]2 \end{array}}$ $Q_{(3.2.1),2n}$\ \ &\  $Q_{(3.2.1),n}$ & $1$\\
$4$ &
$\phantom{\begin{array}{@{}l@{}}1\\[-1.6ex]2 \end{array}}$ $Q_{(3.2.1),2n}$\ \ &\  $Q_{(3.2.2),n}$ & $\frac{13}{29}$\\
$5$ &
$\phantom{\begin{array}{@{}l@{}}1\\[-1.6ex]2 \end{array}}$ $Q_{(3.2.1),2n}$\ \ &\  $Q_{(3.2.3),n}$ & $\frac{1}{3}$\\
$6$ &
$\phantom{\begin{array}{@{}l@{}}1\\[-1.6ex]2 \end{array}}$ $Q_{(3.2.2),2n}$\ \ &\  $Q_{(3.1.1),n}$ & $\frac{168}{235}$\\
$7$ &
$\phantom{\begin{array}{@{}l@{}}1\\[-1.6ex]2 \end{array}}$ $Q_{(3.2.2),2n}$\ \ &\  $Q_{(3.1.3),n}$ & $\frac{28}{15}$\\
$8$ &
$\phantom{\begin{array}{@{}l@{}}1\\[-1.6ex]2 \end{array}}$ $Q_{(3.2.2),2n}$\ \ &\  $Q_{(3.2.2),n}$ & $1$\\
$9$ &
$\phantom{\begin{array}{@{}l@{}}1\\[-1.6ex]2 \end{array}}$ $Q_{(3.2.2),2n}$\ \ &\  $Q_{(3.2.3),n}$ & $\frac{1}{3}$\\
$10$ &
$\phantom{\begin{array}{@{}l@{}}1\\[-1.6ex]2 \end{array}}$ $Q_{(3.2.3),2n}$\ \ &\  $Q_{(3.2.3),n}$ & $1$\\
\hline
\end{tabular}
\end{center}
\caption{Pairs $(Q^{\prime},Q^{\prime\prime})$ of negative
quadrature formulae of order $4$ and the corresponding best
constants $c$, satisfying the assumptions of Theorem \ref{t1}.}
\end{table}

\begin{pf}
All we need is to check that
$\widehat{Q}=(c+1)\,Q^{\prime}-c\,Q^{\prime\prime}\,$ is positive
definite of order $4$. When studying $\,K_4(\widehat{Q};\cdot)\,$ in
the neighborhoods of the endpoints of $[0,1]$, affected by the
formulae for numerical differentiation applied to the construction
of $\,Q^{\prime}\,$ and $\,Q^{\prime\prime}$, we eliminate the
dependence on $n$ by a suitable change of the variable. Away from
these neighborhoods we apply Propositions~\ref{p1} -- \ref{p2} to
obtain a simpler representation of $\,K_4(\widehat{Q};\cdot)\,$. The
verification that $\,K_4(\widehat{Q};\cdot)\,$ does not change its
sign in $(0,1)$ consists of sometimes tedious though elementary
calculations. We therefore decided to present a detailed proof of
only one case, namely, case 9 in Table~1, and point out to some
peculiarities in the other cases.

The interval unaffected by the formulae for numerical
differentiation applied for the construction of $Q^{\prime}$ and
$Q^{\prime\prime}$ in case 9 in Table~1, is $[y_{1,n},y_{n-1,n}]$.
By Proposition~\ref{p2}, for $t\in [y_{1,n},y_{n-1,n}]$ we have
$$
K_4(\widehat{Q};t)=
\frac{c+1}{(2\,n)^4}\,\big[\widetilde{B}_4(2n\,t-\frac{1}{2})-B_4(1/2)\big]
-\frac{c}{n^4}\,\big[\widetilde{B}_4(n\,t-\frac{1}{2})-B_4(1/2)\big]\,.
$$
We shall show that
$$
\varphi(t)=\varphi(c;t):=(c+1)\,\big[\widetilde{B}_4(2n\,t-1/2)-B_4(1/2)\big]
-16c\,\big[\widetilde{B}_4(n\,t-1/2)-B_4(1/2)\big]
$$
is non-negative for every $t\in \mathbb{R}$ if and only if
$c\geq\frac{1}{3}$\,. Since $\,\varphi\,$ is a periodic function
with a period $1/n$, we study its behavior on the interval
$\,[0,1/n]$ only.

Consider first the case $t\in [0,\frac{1}{4n}]\cup
[\frac{3}{4n},\frac{1}{n}]$. If $\,t\in [0,\frac{1}{4n}]$, then we
set $\,t=\frac{1-2u}{4n}\,$, while if $\,t\in
[\frac{3}{4n},\frac{1}{n}]$, then we set $\,t=\frac{3+2u}{4n}\,$,
with $\,u\in [0,\frac{1}{2}]$. In both cases we have
$\,\widetilde{B}_4(2n\,t-1/2)=B_4(u)\,$ and
$\,\widetilde{B}_4(n\,t-1/2)=B_4((2u+1)/4)$, therefore
\[
\begin{split}
\varphi(c;t)&=(c+1)\,\big[B_4(u)-B_4(1/2)\big]
-16c\,\big[B_4((2u+1)/4)-B_4(1/2)\big]\\
&=\frac{(2u-1)^2}{64}\Big[c-\frac{1+4u(1-u)}{6}\Big]\,, \quad u\in
[0,1/2]\,.
\end{split}
\]
The latter expression is non-negative for every $\,u\in[0,1/2]\,$ if
and only if $c\geq 1/3$.

Next, we consider $\,\varphi(t)\,$ with $\,t\in
[\frac{1}{4n},\frac{3}{4n}]$. For $\,t\in
[\frac{1}{4n},\frac{1}{2n}]\,$  we set $\,t=\frac{1-u}{2n}$, while
for $\,t\in [\frac{1}{2n},\frac{3}{4n}]\,$ we set
$\,t=\frac{1+u}{2n}$, with $\,u\in [0,\frac{1}{2}]$. In both cases,
we have $\,\widetilde{B}_4(2n\,t-1/2)=B_4(u+1/2)\,$ and
$\,\widetilde{B}_4(n\,t-1/2)=B_4(u/2)$, therefore
\[
\begin{split}
\varphi(c;t)&=(c+1)\,\big[B_4(u+1/2)-B_4(1/2)\big]
-16c\,\big[B_4(u/2)-B_4(1/2)\big]\\
&=\frac{c}{48}\,\big(8u^3-9u^2+2\big)-\frac{1}{48}\,u^2(1-2u^2)\,,
\quad u\in [0,1/2]\,.
\end{split}
\]
As $\varphi$ is an increasing function of $c$ and $\varphi(1/3;t)=
(3u^4+4u^3-6u^2+1)/72>0$, $u\in [0,1/2]$, we conclude that
$\varphi(t)\geq 0$ in that case, too, provided $c\geq 1/3$.
Consequently, for $c\geq 1/3$ and $t\in [y_{1,n},y_{n-1,n}]$,
$\,K_4(\widehat{Q};t)=(2n)^{-4}\,\varphi(c;t)\geq 0$.

Since $\,\widehat{Q}=(c+1)Q^{\prime}-c\,Q^{\prime\prime}\,$ is a
symmetrical quadrature formula, it remains to verify (with $c=1/3$)
that $\,K_4(\widehat{Q};t)\geq 0\,$ for $\,t\in [0,y_{1,n}]$. As
similar verifications were repeatedly performed in the preceding
section, here we omit the details.

Let us now briefly comment on the other pairs of quadratures in
Table 1. The restriction on $c$ for the pairs of quadratures in
lines $1 - 4$, $6 - 8$ and $10$ of Table~1 comes from the fact that
a closed-type quadrature formula $\,\widehat{Q}\,$ can by positive
definite of order $4$ only if the coefficient of $f(0)$ in
$\,\widehat{Q}\,$ is non-negative, a fact that easily follows from
the explicit form of $\,K_4(\widehat{Q};t)\,$, see \eqref{e3}.
Actually, the values of $c$ in Table~1 in these cases are those, for
which $\widehat{Q}$ is of open type; it turns out that these values
of $c$ secure the positive definiteness of $\,\widehat{Q}\,$.\qed
\end{pf}

\begin{thm}\label{t3} The assumptions of Theorem \ref{t1} are
fulfilled for the pairs $\,(Q^{\prime}, Q^{\prime\prime})\,$ of
positive definite quadrature formulae and with the best possible
constants $c$, given in Table~2.
\end{thm}

\begin{table}[h]\label{tab:2}
\begin{center}
\begin{tabular}{|c|c|c|c|}
\hline No. & $\phantom{\begin{array}{@{}l@{}}1\\[-1ex]2 \end{array}}$\ \ \ \  $Q^{\prime}$\ \ \ \ \ &\ \ \ \ \
$Q^{\prime\prime}$\ \ \ \ &\ \ \ \ $c$\ \ \ \ \\
\hline $1^{\prime}$ &
$\phantom{\begin{array}{@{}l@{}}1\\[-1.6ex]2 \end{array}}$ $Q_{(3.3.1),2n}$\ \ &\  $Q_{(3.3.1),n}$ & $1.104931$\\
$2^{\prime}$ &
$\phantom{\begin{array}{@{}l@{}}1\\[-1.6ex]2 \end{array}}$ $Q_{(3.3.2),2n}$\ \ &\  $Q_{(3.3.1),n}$ & $\frac{1}{3}$\\
$3^{\prime}$ &
$\phantom{\begin{array}{@{}l@{}}1\\[-1.6ex]2 \end{array}}$ $Q_{(3.3.2),2n}$\ \ &\  $Q_{(3.3.2),n}$ & $1.803456$\\
$4^{\prime}$ &
$\phantom{\begin{array}{@{}l@{}}1\\[-1.6ex]2 \end{array}}$ $Q_{(3.3.2),2n}$\ \ &\  $Q_{(3.3.3),n}$ & $1.088270$\\
$5^{\prime}$ &
$\phantom{\begin{array}{@{}l@{}}1\\[-1.6ex]2 \end{array}}$ $Q_{(3.3.2),2n}$\ \ &\  $Q_{(3.4.2),n}$ &
$1.207773$\\
$6^{\prime}$ &
$\phantom{\begin{array}{@{}l@{}}1\\[-1.6ex]2 \end{array}}$ $Q_{(3.3.3),2n}$\ \ &\  $Q_{(3.3.1),n}$ & $\frac{1}{3}$\\
$7^{\prime}$ &
$\phantom{\begin{array}{@{}l@{}}1\\[-1.6ex]2 \end{array}}$ $Q_{(3.3.3),2n}$\ \ &\  $Q_{(3.3.3),n}$ & $1.601589$\\
$8^{\prime}$ &
$\phantom{\begin{array}{@{}l@{}}1\\[-1.6ex]2 \end{array}}$ $Q_{(3.3.3),2n}$\ \ &\  $Q_{(3.4.2),n}$ &
$1.828256$\\
\hline
\end{tabular}
\end{center}
\caption{Pairs $(Q^{\prime},Q^{\prime\prime})$ of positive definite
quadrature formulae of order $4$ and the corresponding best
constants $c$, satisfying the assumptions of Theorem \ref{t1}.}
\end{table}

\begin{pf}
We have to verify that
$\widehat{Q}=(c+1)\,Q^{\prime}-c\,Q^{\prime\prime}\,$ is negative
definite of order $4$. Two kinds of violation of the requirement
$\,K_4(\widehat{Q};t)\leq 0\,$ may occur while decreasing $c$:
\begin{enumerate}[~1)~~]
\item
The requirement is first violated inside the neighborhoods of the
endpoints of $[0,1]$, affected by the formulae for numerical
differentiation applied to the construction of $\,Q^{\prime}\,$ and
$\,Q^{\prime\prime}$. Then the best constant $c$ is a numerically
computed zero of the resultant of a quintic polynomial, with which
the corresponding (re-scaled) Peano kernels coincides.
\item
The requirement is first violated away of these neighborhoods.
There, we exploit Propositions \ref{p3} and \ref{p4} to obtain a
simpler form of $\,K_4(\widehat{Q};\cdot)\,$. Such a situation
occurs in the cases $2^{\prime}$ and $5^{\prime}$ in Table~2.
\end{enumerate}

Here we consider in details only case $2^{\prime}$ in Table 2. The
interval not affected by the formulae for numerical differentiation
applied to the construction of $Q^{\prime}$ and $Q^{\prime\prime}$
is $[y_{1,n},y_{n-1,n}]$. By Proposition \ref{p3}, for $t\in
[y_{1,n},y_{n-1,n}]$ we have
$$
K_4(\widehat{Q};t)=\frac{c+1}{(2n)^4}\,\big[\widetilde{B}_4(2n\,t)-B_4(0)\big]
-\frac{c}{(n)^4}\,\big[\widetilde{B}_4(n\,t)-B_4(0)\big]=:\psi(c;t)\,.
$$
Since $\,\psi(t)=\psi(c;t)\,$ is a periodic function with a period
$1/n$, we may restrict the study of its behavior to the interval
$\,[\frac{1}{2n},\frac{3}{2n}]$.

If $t\in [\frac{1}{2n},\frac{1}{n}]$, we set $\,t=(1/2+u)/n$, $u\in
[0,1/2]$, whence $\widetilde{B}_4(2n\,t)=B_4(2u)$ and
$\widetilde{B}_4(n\,t)=B_4(u+1/2)$. Then
$$
\psi(c;(1/2+u)/n)=
\frac{(1-2u)^{2}}{384\,n^4}\,\big[4u^2-(4u+1)c\big]\,,\qquad u\in
[0,1/2]\,,
$$
and it is non-positive for every $u\in [0,1/2]$ if and only if
$c\geq 1/3$.

If $\,t\in [\frac{1}{n},\frac{3}{2n}]$, we set $\,t=(u+1)/n$,
$\,u\in [0,1/2]$, then $\widetilde{B}_4(2n\,t)=B_4(2u)$ and
$\widetilde{B}_4(n\,t)=B_4(u)$. Now
$$
\psi(c;(u+1)/n)=
\frac{u^2}{96\,n^4}\,\big[(1-2u)^2-(3-4u)c\big]\,,\qquad u\in
[0,1/2]\,,
$$
and it is non-positive for every $u\in [0,1/2]$ if and only if
$c\geq 1/3$\,.

Thus, $\,K_4(\widetilde{Q};t)\leq 0\,$  for $t\in
[y_{1,n},y_{n-1,n}]$ if and only if $c\geq 1/3$. Moreover, $c=1/3$
is the smallest value of $c$ for which
$\,\widehat{Q}=(c+1)\,Q^{\prime}-c\,Q^{\prime\prime}\,$ can be
negative definite of order $4$, where
$(Q^{\prime},Q^{\prime\prime})$ is \emph{any} pair of positive
definite quadrature of order $4$, constructed via the scheme
described in Section 3.3.

Since $\,\widehat{Q}=(c+1)\,Q^{\prime}-c\,Q^{\prime\prime}\,$ is
symmetric, it remains to show, with $c=1/3$, that
$\,K_4(\widehat{Q};t)\leq 0$ for $t\in [0,y_{1,n}]$. The latter is
equivalent to
$$
g(u):=u^4-\frac{1}{3}\,u^3-\frac{64}{27}\,\Big(u-\frac{1}{8}\Big)_{+}^{3}
+2\,\Big(u-\frac{1}{6}\Big)_{+}^{3}
+\frac{8}{9}\,\Big(u-\frac{1}{4}\Big)_{+}^{3}
-\Big(u-\frac{1}{3}\Big)_{+}^{3}\leq 0
$$
for $u\in [0,1/2]$, and it is easily verified to be true. \qed
\end{pf}
%================================================================
\section{Numerical examples}
We have tested the efficiency of the a posteriori error estimates in
Theorem~\ref{t1} for some pairs of quadrature formulae
$(Q^{\prime},Q^{\prime\prime})$ in Tables~1 and 2, with functions
$$
f(x)=e^{x}\,,\qquad
g(x)=-\frac{e^{-x}\,\log\Big(\frac{1+x}{2}\Big)}{\sqrt{1+x}}\,,
$$
which both are $4$-convex, and also have been used in the tests in
\cite{GS:1972}.

In Table 3, the enumeration of the lines corresponds to that in
Tables~1 and 2, and $UEB(Q^{\prime})$ and $UEB(Q^{\prime\prime})$
stand for the upper bounds for $|R[Q^{\prime};\cdot]|$ and
$|R[Q^{\prime\prime};\cdot]|$, provided by Theorem \ref{t1} (ii),
(iii), i.e.,
$$
UEB(Q^{\prime}):=c\,|Q^{\prime}[\cdot]-Q^{\prime\prime}[\cdot]|,\qquad
UEB(Q^{\prime\prime}):=(c+1)\,|Q^{\prime}[\cdot]-Q^{\prime\prime}[\cdot]|\,.
$$
The numerical value of $I[f]$ is $e-1=1.71828182845905\ldots$, which
allows us to evaluate the error overestimation factors
$$
EOF(Q^{\prime}):=\frac{UEB(Q^{\prime})}{\big|e-1-Q^{\prime}[f]\big|},\qquad
EOF(Q^{\prime\prime}):=\frac{UEB(Q^{\prime\prime})}{\big|e-1-Q^{\prime\prime}[f]\big|}.
$$

\FloatBarrier
\begin{table}[!h]\label{tab:3}
\begin{center}
{\small
\begin{tabular}{|c|c|c|c|c|c|c|c|c|}
\hline No.  & function & $n$ & $UEB(Q^{\prime})$ &
$UEB(Q^{\prime\prime})$ & $EOF(Q^{\prime})$ &
$EOF(Q^{\prime\prime})$\\
\hline
\multirow{4}{*}{$4$} & $\phantom{\begin{array}{@{}l@{}}1\\[-1ex]1\end{array}}$
\multirow{2}{*}{$f$} & 16 &  $1.308\times 10^{-8}$ & $4.226\times 10^{-8}$ & 6.813 & 1.359\\
& & 32 & $8.272\times 10^{-10}$ &  $2.672\times 10^{-9}$ & 6.768 & 1.358 \\
& $\phantom{\begin{array}{@{}l@{}}1\\[-1ex]1\end{array}}$
\multirow{2}{*}{$g$} & 16 & $1.369\times 10^{-7}$ & $4.424\times 10^{-7}$ & -- & --\\
& & 32 & $8.749\times 10^{-9}$ & $2.827\times 10^{-8}$ & -- & -- \\
\hline
\multirow{4}{*}{$5$} & $\phantom{\begin{array}{@{}l@{}}1\\[-1ex]1\end{array}}$
\multirow{2}{*}{$f$} & 16 & $9.973\times 10^{-9}$ & $3.989\times 10^{-8}$ & 5.195 & 1.253 \\
& & 32 & $6.228\times 10^{-10}$ & $2.491\times 10^{-9}$ & 5.096 & 1.251 \\
& $\phantom{\begin{array}{@{}l@{}}1\\[-1ex]1\end{array}}$
\multirow{2}{*}{$g$} & 16 & $1.066\times 10^{-7}$ & $4.264\times 10^{-7}$ & -- & -- \\
& & 32 & $6.662\times 10^{-9}$ & $2.665\times 10^{-8}$ & -- & -- \\
\hline \multirow{4}{*}{$9$} & $\phantom{\begin{array}{@{}l@{}}1\\[-1ex]1\end{array}}$
\multirow{2}{*}{$f$} & 16 & $9.957\times 10^{-9}$ & $3.983\times 10^{-8}$ & 5.061 & 1.251 \\
& & 32 & $6.223\times 10^{-10}$ & $2.489\times 10^{-9}$ & 5.030 & 1.250 \\
& $\phantom{\begin{array}{@{}l@{}}1\\[-1ex]1\end{array}}$
\multirow{2}{*}{$g$} & 16 & $1.063\times 10^{-7}$ & $4.251\times 10^{-7}$ & -- & -- \\
& & 32 & $6.652\times 10^{-9}$ & $2.661\times 10^{-8}$ & -- & -- \\
\hline
\multirow{4}{*}{$2^\prime$} &$\phantom{\begin{array}{@{}l@{}}1\\[-1ex]1\end{array}}$
\multirow{2}{*}{$f$} & 16 & $1.128\times 10^{-8}$ & $4.512\times 10^{-8}$ & 5.063 & 1.251 \\
& & 32 & $7.082\times 10^{-10}$ & $2.833\times 10^{-9}$ & 5.031 & 1.250 \\
& $\phantom{\begin{array}{@{}l@{}}1\\[-1ex]1\end{array}}$
\multirow{2}{*}{$g$} & 16 & $1.195\times 10^{-7}$ & $4.780\times 10^{-7}$ & -- & -- \\
& & 32 & $7.539\times 10^{-9}$ & $3.016\times 10^{-8}$ & -- & -- \\
\hline
\multirow{4}{*}{$4^\prime$} & $\phantom{\begin{array}{@{}l@{}}1\\[-1ex]1\end{array}}$
\multirow{2}{*}{$f$} & 16 & $3.596\times 10^{-8}$ & $6.899\times 10^{-8}$ & 16.138 & 1.956 \\
& & 32 & $2.285\times 10^{-9}$ & $4.384\times 10^{-9}$ & 16.232 & 1.957 \\
& $\phantom{\begin{array}{@{}l@{}}1\\[-1ex]1\end{array}}$
\multirow{2}{*}{$g$} & 16 & $3.732\times 10^{-7}$ & $7.162\times 10^{-7}$ & -- & -- \\
& & 32 & $2.406\times 10^{-8}$ & $4.617\times 10^{-8}$ & -- & -- \\
\hline
\multirow{4}{*}{$6^\prime$} & $\phantom{\begin{array}{@{}l@{}}1\\[-1ex]1\end{array}}$
\multirow{2}{*}{$f$} & 16 & $1.128\times 10^{-8}$ & $4.511\times 10^{-8}$ & 5.035 & 1.251 \\
& & 32 & $7.080\times 10^{-10}$ & $2.832\times 10^{-9}$ & 5.017 & 1.250 \\
& $\phantom{\begin{array}{@{}l@{}}1\\[-1ex]1\end{array}}$
\multirow{2}{*}{$g$} & 16 & $1.194\times 10^{-7}$ & $4.777\times 10^{-7}$ & -- & -- \\
& & 32 & $7.537\times 10^{-9}$ & $3.015\times 10^{-8}$ & -- & -- \\
\hline
\end{tabular}}
\end{center}
\caption{Upper bounds for $|R[Q^{\prime};\cdot]|$ and
$|R[Q^{\prime\prime};\cdot]|$, with $(Q^{\prime},Q^{\prime\prime})$
being selected pairs of quadrature formulae from Tables~1 and 2, and
the corresponding error overestimation factors.}
\end{table}
\FloatBarrier

Table~3 depicts the error bounds of six pairs of definite (of the
same kind) quadrature formulae, obtained through Theorem~\ref{t1}.
Although the error bounds provided by the Peano kernel methods may
well overestimate the actual error, we observe here that the error
overestimation factor for the integrand $f$ ranges between $1.250$
and $1.957$ for $Q^{\prime\prime}$, and between $5.017$ and $16.
232$ for $Q^{\prime}$. A conclusion can be drawn also that the error
overestimation factor of $Q^{\prime}$ is greater than the error
overestimation factor of $Q^{\prime\prime}$, although $Q^{\prime}$
provides a better approximation than $Q^{\prime\prime}$. For the
pairs of quadrature formulae $(Q^{\prime},Q^{\prime\prime})$
appearing in Tables~1 and 2 and not included in Table~3, the error
overestimation factor can reach $2.71$ for $Q^{\prime}$ and $32.265$
for $Q^{\prime\prime}$.

Another (and, in fact, frequently used) approach for obtaining error
bounds of definite quadrature formulae is through their error
constants. However, this approach assumes knowledge about the
magnitude of a certain derivative of the integrand, which may not be
available. Here we have $\Vert f^{(4)}\Vert_{C[0,1]}=e$, and hence
an alternative error overestimation factor for a definite quadrature
formula $Q$ of order $4$,
$$
EOF_1(Q):=\frac{e\,|c_4(Q)|}{\big|e-1-Q[f]\big|}.
$$
For the definite quadrature formulae obtained in Section~3, $EOF_1$
varies (rather slightly) between $1.53$ and $1.59$.

So far, we focused on the application of Theorem \ref{t1} for
derivation of error bounds for pairs of definite quadrature formulae
of the same kind. Of course, one should not neglect the classical
approach for obtaining error inclusions through pairs of definite
quadrature formulae of opposite kinds.

As an example, let us consider, e.g., the pair
$(Q^{\prime},Q^{\prime\prime})=(Q_{(3.1.3),n}, Q_{(3.3.3),n})$ of a
negative and a positive definite quadrature formula of order $4$.
$Q^{\prime}$ and $Q^{\prime\prime}$ make use of total $n+7$ nodes.
Following Schmeisser \cite{GS:1972}, we set
$$
\widetilde{I}^{-}:=Q^{\prime},\quad
\widetilde{I}^{+}:=Q^{\prime\prime},\quad
M:=\frac{\widetilde{I}^{-}+\widetilde{I}^{+}}{2},\quad
F:=\frac{\vert\widetilde{I}^{-}-\widetilde{I}^{+}\vert}{2},
$$
thus, for $4$-convex (concave) integrands, $F$ provides an upper
bound for the error of the approximation of the definite integral by
$M$.

\begin{table}[h]\label{tab:4}
\begin{center}
\begin{tabular}{|c|c|c|c|}
\hline function & $\phantom{\begin{array}{@{}l@{}}1\\[-1ex]2 \end{array}}$ $n\ $ &  $M$\ \ &\ \
$F$\ \ \\
\hline \multirow{4}{*}{$f$} &
$\phantom{\begin{array}{@{}l@{}}1\\[-1ex]1\end{array}}$ $12\ \ $& $1.71828183227$ & $1.141\times 10^{-7}$\\
 & $\phantom{\begin{array}{@{}l@{}}1\\[-1ex]1\end{array}}$ $28\ \ $& $1.71828182838$ & $3.732\times 10^{-9}$\\
 & $\phantom{\begin{array}{@{}l@{}}1\\[-1ex]1\end{array}}$ $60\ \ $& $1.71828182845$ & $1.747\times 10^{-10}$\\
\hline \multirow{4}{*}{$g$} &
$\phantom{\begin{array}{@{}l@{}}1\\[-1ex]1\end{array}}$ $12\ \ $& $0.20618061399$ & $1.234\times 10^{-6}$\\
 & $\phantom{\begin{array}{@{}l@{}}1\\[-1ex]1\end{array}}$ $28\ \ $& $0.20618051587$ & $4.050\times 10^{-8}$\\
 & $\phantom{\begin{array}{@{}l@{}}1\\[-1ex]1\end{array}}$ $60\ \ $& $0.20618051540$ & $1.885\times 10^{-9}$\\
\hline
\end{tabular}
\caption{Approximation of $I[f]$ and $I[g]$ by the mean value $M$
and error bounds of the pair $(Q_{(3.1.3),n}, Q_{(3.3.3),n})$ of
definite quadrature formulae of opposite kinds, $n=12,\, 28,\, 60$.}
\end{center}
\end{table}

The values $12,\,28,\, 60$ of $n$ in Table~4 correspond to
$19,\,35,\,67$ nodes used in total by $Q^{\prime}$ and
$Q^{\prime\prime}$, and also to the values $16,\,32,\, 64$ in
\cite[Table 2]{GS:1972}. As is seen, the approximation error $F$
there and in Table~4 behaves similarly.
%======================================================================
\section{Remarks}
\textbf{1.} In \cite{GS:1972} Schmeisser proposed two sequences of
asymptotically optimal positive definite quadrature formulae of
order $4$, which are of open type, i.e., do not involve evaluations
of the integrand at the end-points. It is worth noticing that these
quadrature formulae can be obtained via \eqref{e11} with a slight
modification of $D_1$ (and its reflected variant $\widetilde{D}_1$).
Namely, formulae (45) and (47) in \cite{GS:1972} are obtained
through \eqref{e11} with
$D_1[f]=D_1(x_{0,n},x_{1,n},x_{2,n},x_{3,n},x_{4,n})[f]$ and
$D_1[f]=D_1(x_{0,n},y_{1,n},x_{1,n},x_{2,n},x_{3,n})[f]$,
respectively. Here, $D_1(\mathbf{t})[f]$ stands for the five-point
formula approximating $f^{\prime}(0)$ with nodes
$\mathbf{t}=(0,t_1,t_2,t_3,t_4,t_5)$ and with a \emph{fixed}
coefficient, equal to $-6n$, in front of $f(0)$. With
$$
D_1[f]=-6n\,f(x_{0,n})+\frac{46n}{3}\,f(y_{1,n})-17n\,f(x_{1,n})+
10n\,f(y_{2,n})-\frac{7n}{3}\,f(x_{2,n})
$$
we obtain through \eqref{e11} an $(n+3)$-point ($n\geq 5$) positive
definite quadrature formula of order $4$,
\[
\begin{split}
Q_{n+3}[f]=&\frac{23}{18n}\,[f(y_{1,n})+f(y_{n,n})]
-\frac{5}{12n}\,[f(x_{1,n})+f(x_{n-1,n})]\\
&\!+\!\frac{5}{6n}[f(y_{2,n})\!+\!f(y_{n-1,n})]
\!+\!\frac{29}{36n}[f(x_{2,n})\!+\!f(x_{n-2,n})]
\!+\!\frac{1}{n}\sum_{k=3}^{n-3}f(x_{k,n})
\end{split}
\]
with error constant
$$
c_4(Q_{n+3})=\frac{1}{720\,n^4}\,\Big(1+\frac{55}{4n}\Big)\,.
$$
Compared to the error constants of quadrature formulae (45) and (47)
in \cite{GS:1972} when using the same number of nodes, say, $m$, for
$m\geq 23$ the error constant of the above quadrature formula is
better, i.e., smaller. Yet, it is worse compared to the error
constant of the $m$-point Gaussian quadrature formula $Q_m^{G}$ for
the space of cubic splines with double equidistant knots, which has
been constructed in \cite{GN:1996} and where we have (see
\cite[Corollary 2.3]{GN:1996}, roughly,
$$
c_4(Q_{m}^G)=\frac{1}{720\,(m-1)^4}\,\Big(1-\frac{1.30435}{m-1}\Big)\,.
$$
\smallskip

\textbf{2.} One may wonder why Table~1 does not contain pairs of
negative definite quadrature formulae of order $4$ of the type
$(Q^{\prime},Q^{\prime\prime})=(Q_{(3.1.*),2n}, Q_{(3.1.*),n})$ or
$(Q^{\prime},Q^{\prime\prime})=(Q_{(3.1.*),2n}, Q_{(3.2.*),n})$. The
reason is that, with the above combinations, quadrature formula
$\widehat{Q}=(c+1)Q^{\prime}-c\,Q^{\prime\prime}$ cannot be positive
definite of order $4$ with $c>0$. Indeed, in the first case,
according to Proposition~\ref{p1}, away from the neighborhoods of
the end-points of $[0,1]$, affected by the formulae for numerical
differentiation applied to the construction of $Q^{\prime}$ and
$Q^{\prime\prime}$, we have
$$
K_4(\widehat{Q};t)=\frac{c+1}{(2n)^4}\,\big[\widetilde{B}_4(2n\,t)-B_4(1/2)\big]
-\frac{c}{(n)^4}\,\big[\widetilde{B}_4(n\,t)-B_4(1/2)\big]\,,
$$
and $K_4(\widehat{Q};y_{k,n})<0$ as the first term is negative while
the second term vanishes.

In the second case, by Propositions~\ref{p1} and \ref{p2} we have
away from the end-points
$$
K_4(\widehat{Q};t)=\frac{c+1}{(2n)^4}\,\big[\widetilde{B}_4(2n\,t)-B_4(1/2)\big]
-\frac{c}{(n)^4}\,\big[\widetilde{B}_4(n\,t-1/2)-B_4(1/2)\big]\,,
$$
and $K_4(\widehat{Q};x_{k,n})<0$ as the first term is negative while
the second term vanishes.\smallskip

\textbf{3.} For similar reasons, Table~2 cannot contain pairs of
positive definite quadrature formulae of order $4$ of the type
$(Q^{\prime},Q^{\prime\prime})=(Q_{(3.4.*),2n}, Q_{(3.4.*),n})$ or
$(Q^{\prime},Q^{\prime\prime})=(Q_{(3.4.*),2n}, Q_{(3.3.*),n})$.
Indeed, away from the end-points of $[0,1]$, in the first case one
can see on the basis of Proposition~\ref{p4} that
$K_4(\widehat{Q};y_{k,n})>0$, while in the second case
Propositions~\ref{p3} and \ref{p4} imply
$K_4(\widehat{Q};x_{k,n})>0$.\smallskip

\textbf{4.} Perhaps, the first results on monotonicity of the
remainders of quadrature formulae are due to Newman \cite{DN:1974}.
For conditions for monotonicity of the remainders of quadratures, in
particular of the remainders of compound and Gauss-type quadratures,
in terms of their Peano kernels and the resulting exit criteria, we
refer the reader to \cite{KJF:1993a,KJF:1993,GN:1992, FKN:1998}. The
quadrature formulae constructed here are not of compound type, and
the method applied for proving monotonicity of their remainders by
virtue of Theorem~\ref{t1} (i) is close to that applied in
\cite{KN:1995a}, i.e., relies on the existence of common double
zeros of the shifted Bernoulli monosplines.\smallskip

\textbf{5.} Our choice to construct symmetric quadrature formulae
here is for reasons of simplicity only; otherwise, one can apply
different formulae for numerical differentiation for approximating
the derivatives evaluations at the end-points $0$ and $1$, and thus
obtaining non-symmetric definite quadrature formulae of order~$4$.
Needless to say, the approach proposed here is applicable for the
construction of definite quadrature formulae of higher order.

%\section*{Acknowledgements}
%The authors are partially supported by the Sofia University Research
%Fund through Grant 30/2016.

\end{document}